\RequirePackage{etoolbox}\csdef{input@path}{{style/}{graphics/}}
\documentclass[aos,MSNbibl,nameyear,seceqn,dvips]{arximspdf}
\usepackage{multirow}
\usepackage{graphicx}

\doi{10.1214/13-AOS1193} 
\volume{42}
\issue{2}
\pubyear{2014}
\firstpage{563}
\lastpage{591}

\makeatletter

\newtheorem{thmm}{Theorem}[section]
\newtheorem{cor}{Corollary}[section]
\newtheorem{lem}{Lemma}[section]
\newproclaim{exmp}{Example}[section]
\newproclaim{rmk}{Remark}[section]
\makeatother

\begin{document}
\begin{frontmatter}

\title{Inverse regression for longitudinal data}
\runtitle{Inverse regression}

\begin{aug}
\author[a]{\fnms{Ci-Ren} \snm{Jiang}\thanksref{m1}\ead[label=e1]{cirenjiang@stat.sinica.edu.tw}},
\author[b]{\fnms{Wei} \snm{Yu}\ead[label=e3]{yu.wei@gene.com}}
\and
\author[c]{\fnms{Jane-Ling} \snm{Wang}\corref{}\thanksref{m2}\ead[label=e2]{jlwang.ucdavis@gmail.com}}

\thankstext{m1}{Supported in part by NSC Grant NSC
101-2118-M-001-013-MY2, Taiwan.}
\thankstext{m2}{Supported in part by NSF Grants DMS-04-04630,
DMS-09-06813 and DMS-12-28369.}

\runauthor{C.-R. Jiang, W.~Yu and J.-L. Wang}

\affiliation{Academia Sinica, Genentech Inc and University of
California, Davis}

\address[a]{C.-R. Jiang\\
Institute of Statistical Science\\
Academia Sinica\\
Taipei, 115 \\
Taiwan\\
\printead{e1}}

\address[b]{W. Yu\\
Genentech, Inc\\
1 Dna Way\\
South San Francisco, California 94080\\
USA\\
\printead{e3}}

\address[c]{J.-L. Wang\\
Department of Statistics\\
University of California, Davis\\
Davis, California 95616\\
USA\\
\printead{e2}}
\end{aug}

\received{\smonth{8} \syear{2013}}
\revised{\smonth{11} \syear{2013}}

%
\begin{abstract}
Sliced inverse regression
(Duan and Li [\textit{Ann. Statist.} \textbf{19} (1991) 505--530],
Li [\textit{J. Amer. Statist. Assoc.} \textbf{86} (1991) 316--342])
is an appealing dimension
reduction method for regression models with
multivariate covariates. It has been extended by
Ferr{\'e} and Yao [\textit{Statistics} \textbf{37} (2003) 475--488,
\textit{Statist. Sinica} \textbf{15} (2005) 665--683]
and Hsing and Ren [\textit{Ann. Statist.} \textbf{37} (2009) 726--755]
to functional covariates where
the whole trajectories of random functional covariates are completely
observed. The focus of this paper is to develop sliced inverse
regression for
intermittently and sparsely measured longitudinal covariates.
We develop asymptotic theory for the
new procedure and show, under some regularity conditions, that the
estimated directions attain the optimal rate of convergence.
Simulation studies and data analysis are also provided to demonstrate
the performance of our method.
\end{abstract}

%


\end{frontmatter}

\section{Introduction}\label{sec1}

Dimension reduction methods have played a central role in statistical
modeling with recent interest directed to functional and longitudinal
data. In this paper, we focus on the case when the response is a
univariate variable $Y$, but the covariate $X(\cdot)$ is a stochastic
process that is observed intermittently over a time interval $\mathcal
{I}$, possibly at a few follow-up times. Such observed data are often
termed ``longitudinal'' data in the literature in contrast to
``functional'' data, which are observed densely over a period of time,
so essentially one may assume that the entire process $X$ is observed.

To motivate our approach, we first consider dimension reduction
approaches for a $p$-dimensional multivariate covariate $\mathbf{X}$.
There are essentially two paradigms. The first adopts a model that
reduces the dimensionality of the nonlinear components; this includes
projection pursuit regression [\citet{Friedman81}, \citet{Hall89}] and
additive models [\citet{Stone85}, \citet{Hastie90}]. Although flexible,
both approaches assume a certain additivity structure in the model and
dimension reduction is accomplished during the model fitting stage,
often through an iterative backfitting algorithm.
In the second type of dimension reduction approach, one separates the
dimension reduction stage from the modeling stage, so that model
assumptions are not intertwined with effective dimension reduction.
This approach was pioneered by \citet{Li91} and \citet{Duan91}, who
proposed the inverse regression or sliced
inverse regression (SIR) approach, which assumes that the
information about the response contained in the high dimensional
covariates can be
summarized in a low dimensional subspace.

Specifically,
%
\begin{equation}
Y=f\bigl(\beta_1'\mathbf{X},\ldots,\beta_k'
\mathbf{X},\varepsilon\bigr), \label{sirvec}
\end{equation}
where the $\beta$ are $k$ unknown but
nonrandom vectors with $k < p$, $f$ is an arbitrary unknown link
function on
${R}^{k+1}$, and $\varepsilon$ is a
random error independent of $\mathbf{X}$. An alternative and equivalent
form of
(\ref{sirvec}) is
\[
\mbox{conditional on } \bigl\{\beta_1'\mathbf{X},
\ldots,\beta_k'\mathbf{X}\bigr\}, Y \mbox{ is independent
of } \mathbf{X}. \label{sirvec2}
\]

Hence, under model (\ref{sirvec}), the $k$-dimensional variables
$\{\beta_1'\mathbf{X},\ldots,\beta_k'\mathbf{X}\}$ capture all the
information contained in the original $p$-dimensional variable
$\mathbf{X}$ for predicting $Y$. These models work well for
most practical situations, as the most
interesting features of high dimensional data are usually retrievable from
low-dimensional projections. Note that because the link function
$f$ is unknown, the regression coefficients
$\{\beta_1,\ldots,\beta_k\}$ are not identifiable. However, the
subspace spanned by them is identifiable and is
the ``effective dimension
reduction'' (e.d.r.) space. We also call any direction in
the e.d.r. space an e.d.r. direction. The goal here is to estimate
those directions that span the
e.d.r. space. Many approaches have been proposed to estimate those
e.d.r. directions since Li's pioneering work, including \citet{Cook02}
and approaches that are based on higher order moments [\citet
{Cook91}, \citeauthor{Yin02} (\citeyear{Yin02,YinC:03})]. We focus here on Li's SIR approach, due to its
simplicity and originality.

\citet{Li91} showed that under the design condition,
%
\begin{equation}
E\bigl(b'\mathbf{X}|\beta_1'\mathbf{X},
\ldots,\beta_k'\mathbf{X}\bigr) \mbox{ is linear in }
\beta_1'\mathbf{X},\ldots,\beta_k'
\mathbf{X} \label{sircond1}
\end{equation}
for any direction $ b$ in ${R}^{p}$, the covariance matrix $\operatorname
{cov}[E(\mathbf{X}|Y)]$ is degenerate in any
direction which is $\operatorname{cov}(\mathbf{X})$-orthogonal to the e.d.r.
space spanned by $\{\beta_1,\ldots,\beta_k\}$. Here, $\alpha$ and
$\beta$ are $A$-orthogonal means that satisfy $\alpha'A\beta=0$.
Therefore, the
e.d.r. directions $\{b_j\}_{j=1}^k$ can be located through the
generalized eigen-analysis of $\operatorname{cov}[E(\mathbf{X}|Y)]$ with respect
to $\operatorname{cov}(\mathbf{X})$:
\[
\operatorname{cov}\bigl[E(\mathbf{X}|Y)\bigr]b_j=
\lambda_j\operatorname{cov}(\mathbf{X})b_j.
\]

Since this eigen-analysis only involves $E(\mathbf{X}|Y)$, which marginally
is a one-dimensional nonparametric regression problem as compared to
the original regression $E(Y|\mathbf{X})$, a $p$-dimensional
regression problem, we have accomplished the goal of dimension
reduction through an ``inverse regression.'' 
Once the e.d.r. space is
estimated, standard nonparametric smoothing techniques can then be successfully
applied to the $k$-dimensional covariates
$\beta_1'\mathbf{X},\ldots,\beta_k'\mathbf{X}$, provided that $k$
is much smaller than $p$. The goal of dimension reduction is thus
achieved.

So far, we have briefly discussed traditional
dimension reduction for a multivariate covariate $\mathbf{X}$. We will
next explore this concept for a functional covariate, where $\mathbf{X}$
is replaced by a random function $X(\cdot)$ $ \in$ $L_2(\mathcal{I})$
for an
interval $\mathcal{I} \subset R$. The space $L_2(\mathcal{I})$ is a
collection of Borel measurable real value functions on $\mathcal{I}$,
such that $E(\Vert  X\Vert ^2) =  E(\int_\mathcal{I} |X(t)|^2\,dt) <\infty$.
The modified version of dimension reduction model (\ref{sirvec})
for functional data is
%
\begin{equation}
Y=f\bigl(\langle\beta_1, X \rangle,\ldots,\langle\beta_k,X
\rangle,\varepsilon\bigr), \label{sirfun}
\end{equation}
where 
the $k$ unknown
functions $\beta_1,\ldots,\beta_k$ are in $L_2(\mathcal{I})$, $f$
is an arbitrary unknown function on ${R}^{k+1}$, and $\varepsilon$
is independent of $X(t)$. The notation $\langle u,v \rangle$, for any
$u,v\in L_2(\mathcal{I})$, is
defined as   $\langle u, v \rangle=\int_{\mathcal{I}} u(t) v(t)\,dt$.

\citet{Ferre03} were the first to consider such a dimension reduction
model and to extend SIR to functional data, termed functional SIR. In
a subsequent paper [\citet{Ferre05}], they replaced the
slicing approach for inverse regression by a nonparametric
smoothing method. Further refinements and alternative to functional
SIR have been proposed in \citet{FerrY:07}, \citet{ForzC:07},
\citet{CookFY:10}, and \citet
{ChenHM:11}.
\citet{HsinR:09} provided a different formulation for
the inverse regression method for a scenario where the predictor
$X(t)$ is in a reproducing kernel Hilbert space. 

Extending SIR to functional data is nontrivial, due to the complication
of inverting a covariance operator on $L_2(\mathcal{I})$. An assumption
that is essential for all of these works is that complete trajectories
for a sample of $n$ random functions $X_1, \ldots, X_n$
are fully observed. This assumption, however, is typically
not met in longitudinal studies, as subjects can often only be
measured at discrete and scattered time points, which
may be random and may vary from subject to subject.
Thus, while the observed
longitudinal data originate from underlying smooth random functions,
the observed data have intrinsically different features [\citet{Rice04,
Hall05}]. Longitudinal
data are also often sparsely sampled with very few measurements per subject.

Together with the irregular sampling plan, this poses challenges for
the extension of SIR to longitudinal data.
We will overcome this difficulty by borrowing information across all
subjects in the inverse regression step and by applying smoothing to
estimate the inverse regression function, $E(X(t)| Y)$.
Moreover, our proposed procedure, described in
Section~\ref{sec2.2}, although designed for sparse longitudinal data, can
accommodate more densely sampled
longitudinal data as well. Asymptotic results for the new procedure are
presented in Section~\ref{sec2.3}, where Theorem~\ref{thmm1} implies that the
e.d.r. space can be estimated at a rate that corresponds to that of
one-dimensional smoothing, when the data are sparse. This is the
optimal rate attainable for sparse longitudinal data. We also show that
the parametric
$\sqrt{n}$-rate can be achieved by our method for densely sampled
longitudinal data (or functional data). Thus, our approach not only
resolves the difficulty to adapt SIR for longitudinal data but also
provides a unified platform for functional SIR that can handle multiple
types of sampling frequency for the longitudinal measurements.

The rest of the paper is organized as follows. In the next
section, we state the main approaches, the estimating procedure
and the asymptotic properties. A simulation study and an illustrative
data analysis
are presented in Sections~\ref{sec3} and~\ref{sec4}, respectively, to
demonstrate the effectiveness of the proposed approach. Section~\ref{sec5}
contains concluding remarks. The proofs are
relegated to an \hyperref[app]{Appendix}.

\section{Main approaches and results}\label{sec2}

A similar
condition as (\ref{sircond1}) is needed for functional data:
%
\begin{eqnarray}
\label{fsircond1} &&E\bigl(\langle b, X \rangle|\langle\beta_1, X
\rangle,\ldots,\langle\beta_k, X \rangle\bigr) \mbox{ is linear in }
\nonumber
\\[-8pt]
\\[-8pt]
\nonumber
&&\qquad\langle\beta_1, X \rangle,\ldots,\langle\beta_k, X
\rangle \mbox{ for any direction } b \mbox{ in } L_2(
\mathcal{I}).\label{sircond2}
\end{eqnarray}

Let $\Gamma(s,t)$ and $\Gamma_e(s,t)$ denote the covariance operators
of $X(t)$ and\break  $E(X(t)|Y)$, respectively.
Following similar arguments as those in \citet{Li91}, \citet{Ferre03}
imply that under assumption (\ref{fsircond1}) the operator
$\Gamma_e$ is degenerate in any direction $\Gamma$-orthogonal to the
e.d.r. space. Thus, the basis of the e.d.r. space can be recovered through
the $\Gamma$-orthonormal eigenvectors of $\Gamma_e$, associated with
the $k$ largest eigenvalues:
\[
\label{eqneigen} \Gamma_e\beta_j=\lambda_j
\Gamma\beta_j,
\]
where
$\beta_i'  \Gamma \beta_j=1$, if $i=j$, and 0 otherwise.

Provided that $\Gamma^{-1}$ exists, one could perform a
spectral decomposition of the operator $\Gamma^{-1}\Gamma_e$ (by
requiring $\beta_i'  \Gamma \beta_j=1_{\{i=j\}}$, where $1_{\{\cdot\}
}$ is the indicator function), or equivalently of the operator
$\Gamma^{-1/2}\Gamma_e\Gamma^{-1/2}$ (by requiring the orthogonal
eigenvectors to have norm 1) to locate the e.d.r.
directions. However, such an approach poses difficulties for functional data,
since the compact covariance operator $\Gamma$ is not invertible in the
functional case.
A~practical solution to
regularize the estimate of $\Gamma(s,t)$ or $\Gamma_e (s,t)$ was
proposed by \citeauthor{Ferre03} (\citeyear{Ferre03,Ferre05}). Here, we provide an alternative
approach to define $\Gamma^{-1/2}$ that also illuminates the
identifiability issue of the e.d.r. space.

\subsection{Identifiability of the e.d.r. space}\label{sec2.1}

We have mentioned the identifiability issue of the e.d.r.
directions, that is, the individual vectors $\beta_i$ are not
identifiable. However, the goal of the dimension reduction method
is not about estimating the individual directions
$\beta_1,\ldots,\beta_k$ in model (\ref{sirvec}) [or $\beta_1(t),\ldots
,\beta_k(t)$ in model~(\ref{sirfun})], but
rather estimating the e.d.r. space spanned by
$\{\beta_1,\ldots,\beta_k\}$. 
Since the $\Gamma$-orthonormal eigenfunctons of
$\Gamma_e$ are identifiable, with a little abuse of notation we still
use the notation $\beta_i$ to denote those eigenfunctions and regard
them as the targeted e.d.r. directions.

Another identifiability concern for functional inverse regression
is related to the invertibility of the covariance operator
$\Gamma$. A key issue is how to properly define the unbounded operator
$\Gamma^{-1/2}$ so that the
e.d.r. space can be estimated. Under the assumption that
$E(\Vert X\Vert ^4)<\infty$, $\Gamma(s,t)$ is a
self-adjoint, positive semidefinite, Hilbert--Schmidt operator.
Therefore, there exists an orthonormal basis
$\{\phi_i(t)\}_{i=1}^{\infty}$ in $L_2(\mathcal{I})$ such that
$\Gamma(s,t)$ and $X(t)$ can be represented as
\[
\Gamma(s,t)=\sum_{i=1}^{\infty}
\xi_i\phi_i(s)\phi_i(t),
\]
where $\{\xi_i\}$ are the eigenvalues of $\Gamma$ with
corresponding eigenfunctions $\phi_i(t)$, and $\{\xi_i\}$
satisfies $\xi_1\geq\xi_2\geq\cdots\geq\xi_i \geq\cdots\geq
0$, and
\[
X(t)=\mu(t)+\sum_{i=1}^{\infty}A_i
\phi_i(t),
\]
where $\{A_i\}$ are uncorrelated
random coefficients with $E(A_i)=0, E(A_i^2)=\xi_i$.

To define $\Gamma^{-1/2}$ and $\Gamma^{-1/2}\Gamma_e\Gamma^{-1/2}$, we
consider two scenarios:
\begin{longlist}[1.]
\item[1.] If for some $m$, $\xi_m=0$, the problem is finite
dimensional. We can then use the Moore--Penrose generalized inverse
$\Gamma^{+}(s,t)$ instead of $\Gamma^{-1}$, where
$\Gamma^{+}(s,t)=\sum_{i=1}^{m-1}\frac{1}{\xi_i}\phi_i(s)\phi_i(t)$. Hence,
$\Gamma^{-1/2}(s,t) \triangleq\sum_{i=1}^{m-1}\frac{1}{\sqrt{\xi
_i}}\phi_i(s)\phi_i(t)$.

\item[2.] If there are infinitely many positive eigenvalues, then
$\mathop{\lim}_{i\rightarrow\infty} \xi_i =0$ and the
Moore--Penrose generalized inverse does not exit any more. Thus, we
have to consider another way to define the inverse by restricting
the operator to a smaller domain. When the following condition
is satisfied (see Section~1 in the supplementary material
[\citet{JianYW:13s}] for details):
%
\begin{equation}
\label{eq:id_edr} \sum_{i,j=1}^{\infty}
\frac{E^2\{E(A_i|Y)E(A_j|Y)\}}{\xi_i^2\xi
_j}<\infty,
\end{equation}
a similar augment as in \citet{He03} shows that
$\Gamma^{-1/2}\Gamma_e\Gamma^{-1/2}$ is well defined on the range
space of $\Gamma^{1/2}$, which can be represented as
\[
R_{\Gamma^{1/2}} = \biggl\{f\in L_2(\mathcal{I})\dvtx \sum
_i \xi_i^{-1}\bigl|\langle f,
\phi_i\rangle\bigr|^2 <\infty, f \bot\operatorname{ker}(\Gamma) \biggr
\},
\]
so that for $f \in R_{\Gamma^{1/2}}$,
\[
\Gamma^{-1/2}f \triangleq \sum_{i=1}^{\infty}
\xi_i^{-1/2}\langle f,\phi_i\rangle
\phi_i.
\]
Let $\tilde\Gamma^{-1/2}=\Gamma^{-1/2}|_{R_{\Gamma^{1/2}}}$ denote
such an inverse operator, then the directions we obtain from
$\tilde\Gamma^{-1/2}\Gamma_e\tilde\Gamma^{-1/2}$ are still in the
e.d.r. space, since $R_{\Gamma^{1/2}}$ is a subspace of
$L_2(\mathcal{I})$.
\end{longlist}

\begin{rmk} Condition (\ref{eq:id_edr}) is only a sufficient condition
for\break  $\Gamma^{-1/2}\Gamma_e\Gamma^{-1/2}$ to be well defined. It
originates from another sufficient condition (see Section~1 in the
supplementary material [\citet{JianYW:13s}] for details), that for all $k \geq1$,
%
\begin{equation}
\label{s2:eq2} \sum_i \frac{1}{\xi_i} \biggl(
\sum_j \frac{E\{E(A_i|Y)E(A_j|Y)\}}{(\xi
_i\xi_j)^{1/2}}\langle
\phi_j, \eta_k \rangle \biggr)^2 <\infty,
\end{equation}
which is weaker than (\ref{eq:id_edr}) as can be seen by employing the
Cauchy--Schwarz inequality on the left-hand side of (\ref{s2:eq2}).
However, equation (\ref{s2:eq2}) is not easy to interpret, so we focus
on (\ref{eq:id_edr}) for interpretation.
Simple calculations lead to
\begin{eqnarray*}
&&\sum_{i,j=1}^{\infty}  \frac{E^2\{E(A_i|Y)E(A_j|Y)\}}{\xi_i^2\xi_j}
\\
&&\qquad = \sum_{i,j=1}^{\infty} \frac{1}{\xi_i}
\frac{\operatorname
{var}(E(A_i|Y))\operatorname{var}(E(A_j|Y))}{\xi_i\xi_j} \operatorname {corr}^2\bigl(E(A_i|Y),E(A_j|Y)
\bigr). 
\end{eqnarray*}
Hence, (\ref{eq:id_edr}) is guaranteed if $\operatorname
{var}(E(A_i|Y))/\xi_i$ and the correlations between $E(A_i|Y)$ and
$E(A_j|Y)$ decrease fast enough.
The first requirement is satisfied if the information on $Y$ carried by
$A_i$ decreases fast to zero. The second requirement is also not
stringent since $A_i$ and $A_j$ are uncorrelated principal components.
Below, we provide an example to illustrate (\ref{eq:id_edr}) and when
it does not hold.
\end{rmk}

\begin{exmp}
Let the process $X(t)$ have mean zero and a Karhunan--Lo\'{e}ve expansion,
$X(t) = \sum_i A_i \phi_i(t)$, so $E\{X(t)|Y\} = \sum_i E(A_i|Y) \phi_i(t)$.
\begin{itemize}
\item If $\xi_i \,{=}\, E(A_i^2) \,{=}\, 1/i^2$, $\operatorname{var}(E(A_i|Y))/\xi_i
\,{=}\, 1/i^2$ and $\operatorname{corr}^2\{E(A_i|Y)E(A_j|Y)\} \,{=}\break 1/\{
(i+1)^2(j+1)\}$ for $i\neq j$, then $\sum_{i,j} \frac{E^2\{
E(A_i|y)E(A_j|y)\}}{\xi^2_i\xi_j}$ is finite.
%
\item If $\xi_i \,{=}\, E(A_i^2) \,{=}\, 1/i^2$, $\operatorname{var}(E(A_i|Y))/\xi_i
\,{=}\, 1/i$ and $\operatorname{corr}^2\{E(A_i|Y)E(A_j|Y)\} \,{=}\break  1/\{(i+1)^2(j+1)\}
$ for $i\neq j$, then $\sum_{i,j} \frac{E^2\{E(A_i|Y)E(A_j|Y)\}}{\xi
^2_i\xi_j} \,{=}\, \infty$.
\end{itemize}
\end{exmp}


%

Having resolved the theoretical difficulty with the inverse problem, in
practice, the
estimation of $\Gamma^{-1/2}$ will involve some regularization.
We discuss this in the next subsection.

\subsection{The methodology}\label{sec2.2}

In reality, longitudinal data are
sampled discretely at times $T_{ij}$ from a collection of
trajectories $X_i(t), i=1,\ldots,n$, on a compact interval
$\mathcal{I}$. Following common practice, we assume that the $X_i(t)$
are independent realizations from a
smooth random function $X(t)$ in $L_2(\mathcal{I})$. 
The process $X(t)$ has mean function
$\mu(t)$ and covariance operator $\Gamma(s,t)$ and the scalar
response $Y$ relates to the process $X(t)$ through the
relationship described in (\ref{sirfun}).

Let $X_{ij}= X_i(T_{ij})$ be the $j$th
observation of $X_i$ made at time point $T_{ij}$, where
$i=1,\ldots,n$, and $j=1,\ldots,N_i$. The numbers of
observations $\{N_i\}_{i=1}^{n}$ could be prefixed constants or
i.i.d. random variables sampled from $N$, a discrete random
variable with integer values. For generality, we assume that they
are random variables. The ``observation time points'' $\{T_{ij}\}$
are all in the compact interval $\mathcal{I}$ and assumed to
be i.i.d. realizations of a random variable $T$. In case the $T_{ij}$
are not random, that is, the data are sampled according to a prefixed
schedule, our procedure will still
work, as long as these time points are dense in $\mathcal{I}$.
Using vector notation $\mathbf{T}_i=(T_{i1},\ldots,T_{iN_i})$,
$\mathbf{X}_i=(X_{i1},\ldots,X_{iN_i})$, we can see that
$(\mathbf{T}_i,Y_i,\mathbf{X}_i,N_i)$ are i.i.d., and that the $X_i(t)$
and the $T_{ij}$ are independent of each other even if $X_{ij}$ is
correlated with $T_{ij}$.


To estimate the e.d.r. directions, we adopt the idea of inverse
regression. We first construct the estimators $\hat{\Gamma}$ and $\hat
\Gamma_e$ and then estimate the e.d.r. directions by the eigenfunctions
of $\hat\Gamma_e$ associated with $\hat\Gamma$. The specific steps are:
\begin{enumerate}
\item Estimation of $\Gamma_e$.

For a given time point $t$ and $Y=y$, denote
\[
m(t,y)=E\bigl(X (t) | Y=y\bigr).
\]
We assume that $m(t, y)$ is a smooth function, which can thus be
estimated via a two-dimensional smoothing method applied to the
pooled sample $\{X_{ij}\}$ over $\{T_{ij}, Y_i\}$. While any two-dimensional
smoother can be employed, we use the local linear
regression procedure and derive its asymptotic properties in
Lemma~\ref
{lem1}. Specifically, our objective is to minimize
%
\begin{eqnarray}\label{eq:2dsmooth}
L(t,y)&=&\sum_{i=1}^{n}\sum
_{j=1}^{N_i}\bigl\{X_{ij}-
\alpha_0-\alpha _1(t-T_{ij})-
\alpha_2(y-Y_i)\bigr\}^2
\nonumber
\\[-8pt]
\\[-8pt]
\nonumber
&&\hspace*{28pt}{}\times K_2
\biggl(\frac{T_{ij}-t}{h_t},\frac
{Y_i-y}{h_y}\biggr),
\end{eqnarray}
where $K_{2}(\cdot,\cdot)$ denotes a bivariate kernel function as
defined in (A.3) in the next section, and $h_t$ and $h_y$ are the bandwidths
for $ t$ and $y$, respectively.

%
%

Recall that $\Gamma_e= \operatorname{cov}(E (X(t)| Y))$. Once we have estimated
$m(t, y)$ at all points $t$ and $y$, we can estimate the curve
$E(X(t)|Y_i)$, and hence $\Gamma_e$ by the empirical covariance
function:
\[
\hat \Gamma_e(s,t)=\frac{1}{n}\sum
_{i=1}^{n}\hat m(s,Y_i)\hat
m(t,Y_i)-\frac{1}{n}\sum_{i=1}^{n}
\hat m(s,Y_i)\frac{1}{n}\sum_{i=1}^{n}
\hat m(t,Y_i).
\]

\item Estimation of $\Gamma$.

Noting that $\Gamma=E(X(s)X(t))-E(X(s))E(X(t))$, we need to estimate
the mean
function $\mu(t)= E(X(t))$ and the cross-product $\phi(s, t)=
E(X(s)X(t))$.

The mean function $\mu(t)$ can be estimated by a one-dimensional
local linear smoothing method applied to the pooled data
$\{X_{ij}\}$ over all of the locations $\{T_{ij}\}$, that is, finding
the solution:
\[
(\hat\alpha_0, \hat\alpha_1)=\mathop{\arg
\min}_{(\alpha_0, \alpha_1) \in
{R}^2} \sum_{i=1}^n\sum
_{j=1}^{N_i}\bigl\{X_{ij}-
\alpha_0-\alpha_1(t-T_{ij})\bigr\}
^2K_1\biggl(\frac{T_{ij}-t}{h_\mu}\biggr),
\]
where $K_1(\cdot)$ is a univariate kernel function defined in (A.3)
in the next section. The resulting estimate is $\hat
\mu(t)=\hat\alpha_0$.

To estimate $\phi(s,t)$, a two-dimensional local linear smoother
will be applied to the cross-products $\{X_{ij}X_{ik}\}$. Again,
with
\begin{eqnarray*}
&&(\hat\alpha_0, \hat\alpha_1, \hat\alpha_2)
\\
&&\qquad=\mathop{\operatorname{argmin}}_{(\alpha_0, \alpha_1,\alpha_2) \in
{R}^3}\sum_{i=1}^n
\sum_{j, k=1}^{N_i} \bigl\{X_{ij}X_{ik}-
\alpha_0-\alpha_1(t-T_{ij})-
\alpha_2(s-T_{ik})\bigr\}^2
\\
&&\hspace*{91pt}\qquad\quad{}\times K_1\biggl(\frac{T_{ij}-t}{h_\phi}\biggr)K_1\biggl(
\frac{T_{ik}-s}{h_\phi}\biggr),
\end{eqnarray*}
the estimate is $\hat\phi(s,t)=\hat\alpha_0$.
The estimate for $\Gamma$ is
\[
\hat{\Gamma}(s, t) = \hat{\phi} (s, t) - \hat{\mu}(s) \hat{\mu}(t).
\]


\item Estimation of the e.d.r. directions $\beta_j, j=1, \ldots, k$.

Once we have estimated both $\Gamma_e$ and $\Gamma$, the e.d.r.
directions can be estimated through the eigen-analysis of an estimate
of $\Gamma^{-1/2}\Gamma_e\Gamma^{-1/2}$. Since the eigen-analysis for
an operator can only be performed in reality on a discrete time grid,
the actual implementation involves the following steps:

\begin{enumerate}[(a)]

\item[(a)] Discretize $\hat\Gamma_e$ and $\hat\Gamma$ on an
equally-spaced grid $\{t_1, \ldots, t_p\}$ to obtain the $p\times p$
matrices, $\hat\Gamma_{e,p}$
and $\hat\Gamma_{n,p}$; %
\item[(b)] Perform the singular value decomposition on
$\hat\Gamma_{n,p}$
to compute $\hat\Gamma_{n,p}^{-1/2}$; %
\item[(c)] Perform the eigenvalue decomposition of
$\hat\Gamma_{n,p}^{-1/2}  \hat\Gamma_{e,p}  \hat\Gamma_{n,p}^{-1/2}$ to obtain the first $k$ eigenfunctions
$\hat\eta_1,\ldots,\hat\eta_k$ corresponding to the $k$
largest eigenvalues $\hat\lambda_1,\ldots,\hat\lambda_k$;
\item[(d)]
$\hat\beta_j=\hat\Gamma_{n,p}^{-1/2}\hat\eta_j, j=1,\ldots,k$.
\end{enumerate}
\end{enumerate}
Here, $\hat\eta_j$ are the estimates for the standardized
e.d.r. directions $\eta_j$, and $\hat\beta_j$ for $\beta_j$.

\begin{rmk} If measurements are taken according to the same time
schedule for all subjects, it is natural to use this time schedule as
the grid in step 3(a) above. A guiding principle to select the grid
size $p$ for general longitudinal applications is to
choose it large enough to reveal the characteristics of the function,
but not so large as to generate an excessive computational burden or a
large $p$ problem. 
Our experience has been that there is no significance impact on the
performance of the approach unless the choice of $p$ is grossly wrong.

The discretization step in 3(a) provides some level of regularization
for the inverse operator, $\hat\Gamma_{n,p}^{-1/2}$, in step 3(b). We
found that additional regularization through truncating the smallest
eigencomponents of
$\hat\Gamma_{n,p}$ is often helpful. A~rule of thumb that has worked
well in numerical studies is to retain only the first $L$
eigencomponents that explain a desirable fraction of the variation of
$\hat\Gamma_{n,p}$. The strategy we adopt is to start with a large
fraction, say 0.99, and then decrease it gradually until the global
pattern of the estimated direction emerges. Such a scheme also
automatically excludes components of $\hat{\Gamma}_{n, p}$ that have
negative eigenvalues, so the final covariance estimate will be positive
semidefinite and a square root can be taken using the Moore--Penrose
generalized inverse. This approach to reconstruct the covariance
estimate provides the optimal projection onto the space of positive
semi-definite covariance operator as shown in \citet{HallMY:08}.
\end{rmk}

\subsection{Asymptotic properties}\label{sec2.3}

We present the consistency and rates of convergence of the estimated
covariance operators and the e.d.r. directions in this section.
The convergence of the two-dimensional local linear estimator of
$E(X (t)| Y=y)$ and the convergence of the estimated covariance
operators, $\hat\Gamma_e$, on bounded intervals
are key results and of independent interest (Lemmas \ref{lem1} and \ref{lem2}).
The main result on the convergence of
the e.d.r. directions is presented in Theorem~\ref{thmm1}.
All proofs are in the \hyperref[app]{Appendix}, except for the proof of
Lemma~\ref
{lem2}, which is provided in the supplementary material [\citet{JianYW:13s}].

The estimators $\hat\Gamma$ and $\hat\Gamma_e$ have been
constructed by the local linear smoothing method. Therefore, it is
natural to make the standard smoothness assumptions on the second
derivatives of $\Gamma$ and $\Gamma_e$. Assumed that the data
$(\mathbf{T}_i,Y_i,\mathbf{X}_i), i=1,\ldots, n$, have the same
distribution, where
$\mathbf{T}_i=(T_{i1},\ldots,T_{iN_i})$ and
$\mathbf{X}_i=(X_{i1},\ldots,X_{iN_i})$. Notice that
$(T_{ij},Y_i)$ and $(T_{ik},Y_i)$ are dependent but identically
distributed, we assume they have the marginal density
$g(t,y)$. The assumption~(\ref{fsircond1}) and condition (\ref
{eq:id_edr}) are assumed to hold throughout the paper. Additional
assumptions are listed in
(A.1)--(A.8) below.

(A.1) The numbers of observations
$\{N_i\}$ are independent random variables, with
$\stackrel{\mathrm{i.i.d.}}{\sim} N$, where $N$ is a positive integer random variable
with $P(N>1)>0$. Here, we will view $N$ as a random function of
sample size $n$, and $N(n)$ may go to $\infty$ as
$n\rightarrow\infty$. Furthermore, we assume that
$\mathop{\sup}_{n\rightarrow\infty}\frac{E(N(n)^2)}{(EN(n))^2}<\infty$, and
$\mathop{\sup}_{n\rightarrow\infty}\frac
{E(N(n)^4)}{(EN(n)^3EN(n))}<\infty$.\vspace*{1pt}

These conditions are automatically satisfied when
$N_i$'s are uniformly bounded, that is, $N(n) $ is uniformly bounded by a
given positive integer $M$ such that $P(N(n)<M)=1$, $\forall n$.
Therefore, assumption (A.1) is intended for nonsparse data only.


(A.2) As mentioned before, $\{\mathbf{T}_i,Y_i,\mathbf{X}_i,
N_i\}$ are independent. Furthermore,
$\mathbf{T}_i,Y_i,\mathbf{X}_i$ are independent of $N_i$.

(A.3) Let $K_2(\cdot,\cdot)$ be the bivariate kernel function,
which is compactly supported, symmetric and H\"{o}lder
continuous. We further assume that it is a kernel of order
$(|v|,|\kappa|)$, that is,
\[
\sum_{\ell_1+\ell_2 = \ell} \int\int u^{\ell_1}v^{\ell_2}K_2(u,v)
\,du \,dv =
\cases{ 0, &\quad  $0\leq\ell<|\kappa|, \ell\neq|v|,$\vspace*{2pt}
\cr
(-1)^{|v|}|v|!, &\quad $\ell= |v|,$\vspace*{2pt}
\cr
\neq0, & \quad $\ell=|
\kappa|.$}
\]
Similarly, the univariate kernel function
$K_1$ can be defined. We say that $K_1$ is of order $(m,k)$, if
\[
\int u^{\ell}K_1(u)=\cases{ 0, &\quad $0\leq\ell<k, \ell\neq m$,
\vspace*{2pt}
\cr
(-1)^{m}m!, & \quad$\ell=m$,\vspace*{2pt}
\cr
\neq0, & \quad$
\ell=k$.}
\]
Both $K_1$ and $K_2$ are square integrable, that is, $\int
K_i^2\,du<\infty$, $i=1,2$.

In our application, we set
$|v|=0, |\kappa|=2$ for $K_2$ and $m=0, k=2$ for $K_1$, but other order
could be used by properly adjusting the results.


(A.4) Without loss of generality, we assume that $h_t$ and $h_y$
have the same order:
\[
h_t\sim O(h),\qquad h_y\sim O(h).
\]

%


(A.5) The bandwidth $h$ satisfies $\lim_{n\rightarrow\infty} h = 0$,
$\lim_{n\rightarrow\infty} hEN(n) < \infty$, $\lim_{n\rightarrow\infty
} n h^2EN(n) = \infty$, and $\lim_{n\rightarrow\infty} n
h^6EN(n)<\infty$.


(A.6) Let $g(t,y)$ be the density function of $(T_{ij},Y_i)$ and
$g_3(s,t,y)$ be the joint density function of $(T_{ij}, T_{ik}, Y_i)$.
Let $\Psi(s,t,y)$ be the covariance of $X(S)$ and $X(T)$ given $S=s$,
$T=t$ and $Y=y$. We assume that
$g$, $g_3$ and $\Psi$ have continuous and bounded second derivatives
and that $g$ is bounded away from zero.

(A.7) The bandwidth $h_\mu$ satisfies $\lim_{n\rightarrow\infty} h_\mu
= 0$, $\lim_{n\rightarrow\infty} h_\mu EN(n) < \infty$, $ \lim_{n\rightarrow\infty} n h_\mu^2EN(n)= \infty$, and $\lim_{n\rightarrow
\infty} n h_\mu^6EN(n)<\infty$. The bandwidth $h_\phi$ satisfies $\lim_{n\rightarrow\infty} h_\phi= 0$, $\lim_{n\rightarrow\infty} h_\phi
EN(n)^3 < \infty$, $\lim_{n\rightarrow\infty} n h_\phi^2EN(n)^2 =
\infty$, and $\lim_{n\rightarrow\infty} n h_\phi^6EN(n)^2<\infty$.


(A.8) $E(\|X\|^4)<\infty$.


\begin{lem}\label{lem1}
Under assumptions \textup{(A.1)--(A.6)}, we have
%
\begin{eqnarray}\label{eq:bias}\quad
&& E\bigl(\hat{m}(t,y)-m(t,y)|\mathbf{T}_i,Y_i,i =1,
\ldots,n\bigr)
\nonumber
\\[-8pt]
\\[-8pt]
\nonumber
&&\qquad= \frac{\sigma^2}{2}\operatorname{tr}\bigl(H\cdot \mathcal{H}_m(t,y)
\bigr)+o_p\bigl(\operatorname{tr}(H)\bigr),
\\
&&\label{eq:var}
\operatorname{var}\bigl(\hat{m}  (t,y)|\mathbf{T}_i,Y_i,i=1,
\ldots,n\bigr)
\nonumber
\\[-8pt]
\\[-8pt]
\nonumber
&&\qquad = \frac{\Psi(t,t,y)}{nEN(n)|H|^{1/2}}\bigl\{ R(K_2)g(t,y)+\delta
g_3(t,t,y)\bigr\}/g^2(t,y)+o_p(1),
\end{eqnarray}
where $\sigma^2 = \int v^2K_2(v)\,dv$, $H=\operatorname{diag}\{h_t^2,h_y^2\}$, $\mathcal
{H}_m(t,y)$ is the
second derivative of $m(t,y)$, $R(K_2)=\int
K_2^2(u)\,du$, and $\frac{E[N(n)\{N(n)-1\}]}{E\{N(h)\}
h_y}|H|^{1/2}\rightarrow\delta$.

When $EN(n)<\infty$, which is the case for longitudinal data, $\delta$
is zero. Thus, the variance (\ref{eq:var}) can be simplified to
\[
\operatorname{var}\bigl(\hat{m} (t,y)|\mathbf{T}_i,Y_i,i=1,
\ldots,n\bigr) = \frac
{\Psi(t,t,y)}{nEN(n)|H|^{1/2}}R(K_2)/g(t,y)+o_p(1).
\]
\end{lem}

After estimating $E(X(t)|Y)$, $\Gamma_e$ can be estimated empirically.
Specifically,
\[
\hat\Gamma_e(s,t)=\frac{1}{n}\sum
_{i=1}^{n}\Biggl[\hat m(s,Y_i)-
\frac{1}{n}\sum_{j=1}^{n}\hat
m(s,Y_j)\Biggr] \Biggl[\hat m(t,Y_i)-\frac{1}{n}
\sum_{j=1}^{n}\hat m(t,Y_j)
\Biggr].
\]

From Lemma~\ref{lem1}, we obtain the following.

\begin{lem}\label{lem2}
Under the
assumptions of Lemma~\ref{lem1},
\[
\bigl\Vert \hat{\Gamma}_e(s,t)-\Gamma_e(s,t)\bigr\Vert =
O_p \biggl(\frac{1}{\sqrt{n
h^2EN(n)}}+h^2 \biggr).
\]

Here, $\Vert \cdot\Vert $ denotes the operator norm in $L_2$.
\end{lem}

Lemmas \ref{lem1} and~\ref{lem2} implies that we have the same rate of
convergence as the conventional case for smoothing two-dimensional
independent data. Thus, the within subject dependency causes
technical difficulties but one does not pay a price in the convergence
rate.

We also need the convergence of $\hat\Gamma$.

\begin{lem}\label{lem3}
Under assumptions \textup{(A.1)--(A.3)}, and
\textup{(A.6)--(A.8)},
\[
\bigl\Vert \hat{\Gamma}(s,t)-\Gamma(s,t)\bigr\Vert = O_p \biggl(
\frac{1}{\sqrt{nh_\phi^2 EN(n)} }+\frac{1}{\sqrt{nh_\mu
EN(n)}}+(h_\mu+h_\phi)^2
\biggr).
\]
\end{lem}

An immediate application of Lemmas~\ref{lem1}--\ref{lem3} leads
to the following.

\begin{lem} \label{lem4}
Under assumptions \textup{(A.1)--(A.8)},
\begin{eqnarray*}
\bigl\|\hat\Gamma^{-1/2}\hat\Gamma_e\hat\Gamma^{-1/2}-
\Gamma^{-1/2}\Gamma _e\Gamma^{-1/2}\bigr\|&=&
O_p\bigl(\|\hat{\Gamma}-\Gamma\|+\|\hat{\Gamma}_e-
\Gamma_e\|\bigr)
\\
&=&O_p \biggl(\frac{1}{\sqrt{nh_T^2 EN(n)}}+h_T^2
\biggr),
\end{eqnarray*}
where $h_T=h+h_\mu+h_\phi$.
\end{lem}

From Lemma~\ref{lem4}, we can see that the optimal convergence rate is
achieved when the bandwidth $h_T\sim(nEN(n))^{-1/6}$. This is
also the optimal rate for the two-dimensional smoothing step involved in
both $\Gamma$ and $\Gamma_e$.

Once estimates $\hat\Gamma_e$ and $\hat\Gamma$ for $\Gamma_e$ and
$\Gamma$ have been obtained,
we proceed to estimate the
$j$th eigenfunction $\eta_j$ based on:
$\hat\Gamma^{-1/2}\hat\Gamma_e\hat\Gamma^{-1/2}\hat\eta_j=\hat\lambda
_j\hat\eta_j$.
From the perturbation theory for linear operators [\citet{Kato:66},
Chapter VIII], we readily obtain the following.

\begin{cor}\label{cor1}
Under assumptions \textup{(A.1)--(A.8)} and
the assumption that the nonzero eigenvalues $\{\hat\lambda_j\}$
are distinct, the eigenvector $\hat\eta_j$ satisfies
\[
\|\hat\eta_j-\eta_j\|= O_p \biggl(
\frac{1}{\sqrt{nh_T^2 EN(n)}}+h_T^2 \biggr)
\]
for $1\leq j \leq k$.
\end{cor}

The rate of convergence in Corollary~\ref{cor1} is
the same as the rate of convergence for the covariance estimates of
$\Gamma$ and
$\Gamma_e$ stated in Lemmas \ref{lem2} and~\ref{lem3}. These
rates correspond to
the traditional optimal
rate for a two-dimensional smoother, and is a consequence of applying
perturbation theory. While such a rate is optimal to estimate the covariance
operator,
it is not optimal for the estimation of
eigenfunctions, which should be estimable at the optimal rate
for a one-dimensional smoother. With extra technical work, this is
indeed achievable and the result is presented in the next theorem. 

\begin{thmm}\label{thmm1}
Under assumptions \textup{(A.1)--(A.8)}, 
we have
%
\begin{equation}
\Vert \hat\eta_j-\eta_j\Vert = O_p
\biggl(\frac{1}{\sqrt{nhEN(n)}}+h^2+\frac{1}{\sqrt{nh_\phi
EN(n)}}+h_\phi^2
\biggr) \label{etarate}
\end{equation}
for $1\leq j \leq k$. Therefore, three types of optimal convergent
rates emerge:
\begin{longlist}[(iii)]
\item[(i)] when $EN(n)<\infty$ (sparse longitudinal data), the optimal
rate of $\|\hat\eta_j-\eta_j\|$ is
$O_p (n^{-2/5} )$;
\item[(ii)] when $EN(n)h\rightarrow0$, but $EN(n)\rightarrow\infty$,
the optimal rate of $\|\hat\eta_j-\eta_j\|$ is $O_p (n^{-r} )$,
where $2/5<r<1/2$;
\item[(iii)] when $EN(n)h<\infty$ (dense longitudinal data or
functional data), the optimal rate of $\|\hat\eta_j-\eta_j\|$ is
$O_p (n^{-1/2} )$.
\end{longlist}
\end{thmm}

Note that $h$ and $h_\phi$ are the bandwidths used in the
two-dimensional local linear smoothing on $\{X_{ij}(T_{ij},Y_i)\}$ and
$\{X_{ij}X_{ik}\}$, where the optimal order of $h$ and $h_\phi$ is
$(nEN(n))^{-1/6}$ according to Lemmas \ref{lem1}--\ref{lem3}.
Therefore, (\ref{etarate}) demonstrates the need of undersmoothing in
order to estimate $\eta_j$ optimally.

\section{Simulation studies}\label{sec3}

Since our method is applicable to both sparse longitudinal data
and functional data, we evaluate its finite sample performance for both
types of data through simulations. Without loss of generality, we set
the domain
interval $\mathcal{I}$ as $[0,1]$.
Let $X(t)$ be a standard Brownian motion on $[0,1]$; we consider
the following model:
\[
Y=3+\exp\bigl(\bigl\langle\beta(t),X(t)\bigr\rangle\bigr)+\varepsilon,
\]
%
where $\beta(t)=\sqrt{2}\sin(3\pi t/2)$,
and the random error $\varepsilon\sim N(0,0.1^2)$. The standard deviation
of $\varepsilon$ may look small, but the range of $\exp(\langle\beta
(t),X(t)\rangle)$ is around $(0.5,1.5)$, so the signal-noise ratio is
about $10$.

In each run, $n$ sample trajectories, $\{ X_i(t), i=1,\ldots,n\}$, are
generated from
Brownian motion on [0, 1]. This forms the complete data, but for
practical implementation we discretized the data to equally spaced
31 time-points, $\{t_0,t_1,\ldots,t_{30}\}$, with $t_0=0$ and
$t_{30}=1$. Therefore, the actual dense data set is
$\{X_{ij}=X_i(t_{ij}), i=1,\ldots, n, j=1,\ldots, 30\}$ along with
its response $Y_i$. To generate the sparse longitudinal data, we
randomly selected 2 to 10 observations from $\{t_1,t_2,\ldots,t_{30}\}$.
This results in the longitudinal data $(X_{i1},\ldots,X_{iN_i})$ for
the $i$th subject at time points
$(t_{i1},\ldots,t_{iN_i})$, where $N_i$ follows a uniform distribution
on $\{2,3,\ldots,10\}$. The simulation consists of 100 runs and Table~\ref{tablenumresult} summarizes the numerical findings when $n$ is 100
and 200. As a comparison, we also include the results of the smoothed
functional inverse regression approach in \citet{Ferre05}, which is for
complete data.

%
\begin{table}
\caption{Simulation comparison of FY [Ferr{\'e} and Yao (\citeyear{Ferre05})] for complete
data and our procedures for both complete and sparse data. The
comparison is based on the averages of correlations, ISB, IVAR and IMSE
in 100 simulation runs}\label{tablenumresult}
\begin{tabular*}{\textwidth}{@{\extracolsep{\fill}}lccccc@{}}
\hline
\multicolumn{1}{@{}l}{$\bolds{n}$} & \textbf{Data type} & \textbf{Correlation} & \textbf{ISB} &
\textbf{IVAR} & \textbf{IMSE} \\
\hline
{100} & FY (Complete) & 0.7159 & 0.0114 & 0.2008 &
0.2123\\
& Complete & 0.9912 & 0.0043 & 0.0084 & 0.0127\\
& Sparse & 0.8831 & 0.0583 & 0.2823 & 0.3406 \\[3pt]
{200} & FY (Complete) & 0.8218 & 0.0024 & 0.0837 &
0.0861\\
& Complete & 0.9921 & 0.0024 & 0.0092 & 0.0116\\
& Sparse & 0.9438 & 0.0274 & 0.1602 & 0.1876 \\
\hline
\end{tabular*}
\end{table}


The first comparison is based on the correlation between $\langle\beta
(t), X(t)\rangle$ and $\langle\hat\beta(t), X(t)\rangle$, that is, the
correlation between the projection
of $X(t)$ on the real e.d.r. direction and that on the estimated e.d.r.
direction. Averages of those correlations are reported in the third
column of Table~\ref{tablenumresult}. The results suggest that our
approach generally produces high correlations and for complete data
these are larger than those reported in \citet{Ferre05}. 
The remaining comparisons are based on the Integrated Squared
Bias (ISB), Integrated Variance (IVAR) and Integrated Mean Square
Error (IMSE) [or Mean of Integrated Square Error (MISE)].\vadjust{\goodbreak} The \hyperref[app]{Appendix}
contains details of those definitions. The averages of these statistics
over the 100 simulation runs are reported in columns 4--6 of Table~\ref{tablenumresult}. As expected, the results for complete data are better
than those for sparse data and the results for larger sample sizes are
better. For complete data, our procedure generally led to smaller ISB,
IVAR and IMSE than Ferr\'{e} and Yao's.

%
\begin{figure}[b]

\includegraphics{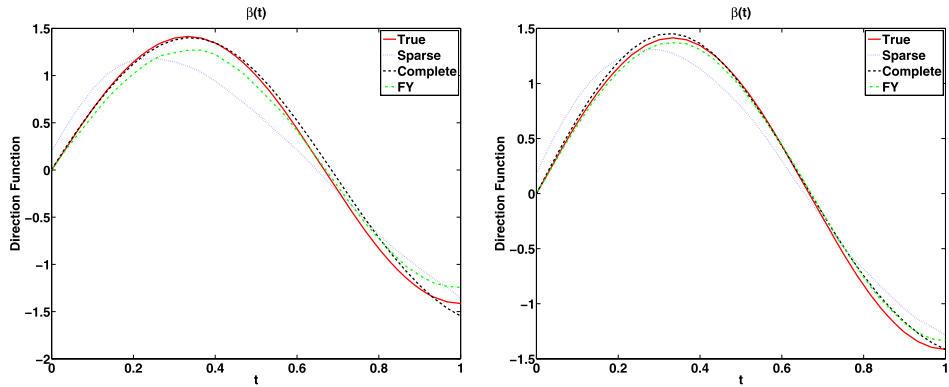}

\caption{{Simulation comparison of the average estimates of $\beta(t)$
for the three methods in Table~\protect\ref{tablenumresult}. The left panel shows the average of
$n$ estimate ($\hat\beta(t)$) for various methods
vs. the target ($\beta(t)$) for $n=100$, and the right panel for $n=
200$.}}\label{simfig}
\end{figure}


In addition to the above global measures, we plot in Figure~\ref{simfig} the mean function for each of the three $\beta$-estimates. The
left panel of Figure~\ref{simfig} shows the average of $\hat\beta
(t)$-functions [dashed line for complete, dotted line for sparse data
and dash-dot line for \citet{Ferre05}] when $n$ = 100 along with the
true $\beta(t)$ (solid
line), the right panel provides the same plot for $n$ = 200. Figure~\ref{simfig} indicates that bias for our approach is comparable to that
reported in \citet{Ferre05} when data are observed completely.
The bias of our approach is significantly reduced for sparse data when
the sample size increases to 200, due to improved estimation of $\Gamma
$ and~$\Gamma_e$.



Upon the request of a referee, we conducted additional simulations with
different sample sizes to check the empirical convergence rate of the
standardized e.d.r. directions ($\eta_k$) through integrated variance
(IVAR). Using the same bandwidths and $N_i=6$ for all $i$, the ratios
of $\sqrt{\operatorname{IVAR}}$ for two consecutive samples (100 vs. 200 or 200
vs. 400) are close to $\sqrt{2}$, which is the square root of the ratio
of sample sizes (see the supplementary material
[\citet{JianYW:13s}]).\vspace*{-2pt}


\section{Data analysis}\label{sec4}

The data set contains the record of the
lifetimes and daily reproduction of female Medflies, the latter quantified
by the number of eggs laid daily for 1000 female
Mediterranean fruit flies. Details about the experimental background
can be
found in \citet{Carey98}. Our goal is to explore the relationship
between the early pattern of fecundity, quantified by the number of
eggs laid per day until day 20, and mortality for each individual fly.
For this reason, we exclude flies that died by day 30 and flies
that did not lay any eggs. The remaining 647 flies have an average
lifetime ($Y$) of 43.9 days with a standard deviation of 11.9 days.
It is assumed that there is an underlying stochastic predicting process
$X(t)$ which quantifies the reproduction pattern and can be
characterized as a fecundity
curve that is sampled through the daily egg counts. The numbers of eggs laid
in the first 20 days are discrete observations of the function
$X(t)$. The objective of our analysis is to find the
e.d.r. directions such that the projection of the fecundity curves
onto the resulting e.d.r. space will carry the key information for
longevity in the regression $E(Y|X)$.

To test the efficiency of our method and to check the effect of sparse
data, we first use the complete information of all 20 days as \textit
{complete/dense data}; and then randomly pick $N_i$ points from each
fly as our \textit{sparse data}, where $N_i$ is uniformly distributed
in $\{2,\ldots, 10\}$. We also applied the approach in \citet{Ferre05}
to the complete data as a comparison. 

\begin{figure}[b]\vspace*{-2pt}

\includegraphics{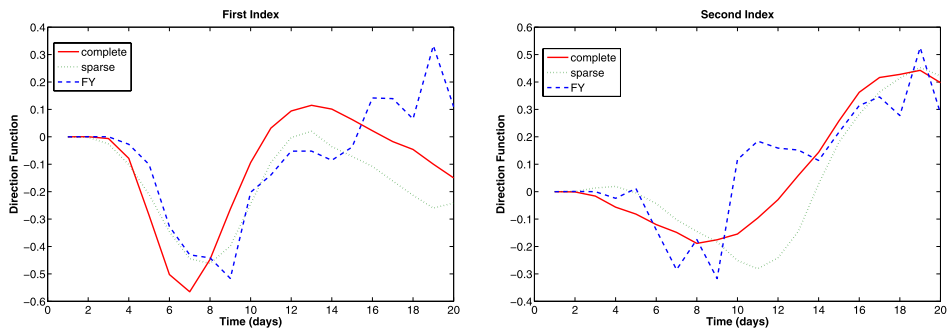}

\caption{Estimated $\beta_1(t)$ and $\beta_2(t)$ from complete
and sparse fecundity data.}
\label{figeggdir1}
\end{figure}


Figure~\ref{figeggdir1} displays the directions estimated by our
approach for both complete and sparse data\vadjust{\goodbreak} and by the method in \citet
{Ferre05} for the complete data only. The directions estimated by \citet
{Ferre05} are less smooth because $\Gamma(s,t)$ was estimated
empirically without smoothing. However, the general trends of these
directions are similar to ours except for the first index after 15 days
($t>15$). The global patterns of the direction estimates by our
approach are similar between the two types of data and the difference
might be due to the difference in the selected bandwidths (a~larger
bandwidth is used for sparse data to compensate for the sparsity, and
this leads to smoother directions).
The estimated $\beta_1(t)$ indicates that daily reproduction during the
period day 4 to day 10 plays an important role in mortality, while the
estimated $\beta_2(t)$ shows the effect of daily reproduction from day
10 to day 20.\vadjust{\goodbreak}

\begin{figure}[b]\vspace*{3pt}

\includegraphics{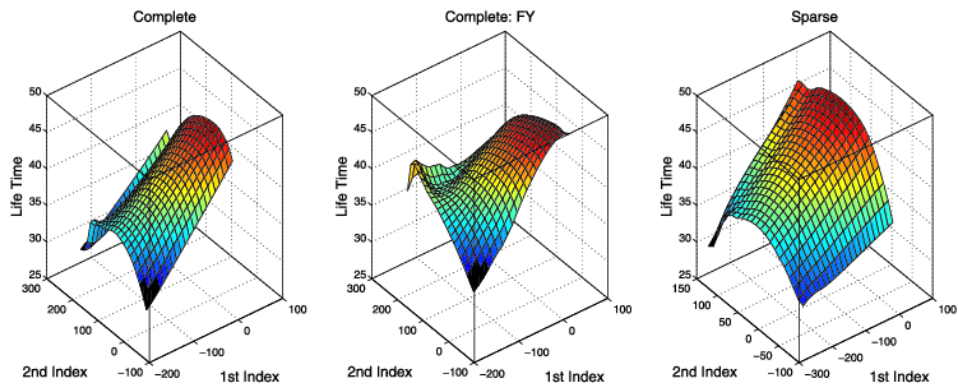}

\caption{Estimated link functions for the fecundity data: the left
panel shows the regression surface when the directions were estimated
by the approach in Ferr{\'e} and Yao (\citeyear{Ferre05}) with complete data; the middle panel
shows the regression surface when the directions were estimated by our
approach with complete data; the right shows the regression surface
when the directions were estimated by our approach with sparse data,
but the indices were calculated from the true complete covariate.}
\label{figegglk}
\end{figure}

%


Since the first two eigenfunctions explain over $90\%$ of the variation
for both sparse and complete data, two directions suffice to summarize
the information contained in the fecundity data to infer lifetime.
We further explore the relation of lifetimes with these two directions
by assuming that the error $\varepsilon$ in model (\ref{sirfun}) is
additive but the regression relation
is unknown. This unknown bivariate regression function is estimated by a
bivariate local linear smoother on the estimated bivariate indices
($\langle\hat\beta_1(t), X_i(t)\rangle$ and $\langle\hat\beta_2(t),
X_i(t)\rangle$). Details regarding the bivariate local linear smoother
are provided in the data analysis section of the supplementary material [\citet{JianYW:13s}].
The estimated regression (link) surfaces are provided in Figure~\ref{figegglk}, where the indices on the right panel were obtained by using
the directions estimated from sparse data but the complete covariate
$X(t)$ was used to calculate the indices. This facilitates a comparison
on the same platform with the other two plots, where the indices were
estimated based on complete data.

Since the estimated e.d.r. directions are not identical, the ranges of
the resulting indices are slightly different. For both sparse and
complete data, lifetimes tend to increase with increasing size of the
first index when the second index is held fixed. 
Lifetime is generally longer when the first index is larger and the
second index is close to its average value. Averages of the square
fitted errors are provided in Table~\ref{tb:fiterr} and are similar for
all three methods, but interestingly our approach for sparse data
performed slightly better than \citeauthor{Ferre05}'s (\citeyear{Ferre05}) approach based on
complete data.

\begin{table}
\caption{Average of the square fitted errors, $\frac{1}{n}\sum_{i=1}^n
(y_i-\hat{y}_i)^2$, of the fecundity data, based on three different
methods: Ferr{\'e} and Yao (\citeyear{Ferre05}), our method with complete data, and our method
with sparse data}\label{tb:fiterr} 
%
\begin{tabular*}{\textwidth}{@{\extracolsep{\fill}}lccc@{}}
\hline
\textbf{Method} & \textbf{Complete: FY} & \textbf{Complete} & \textbf{Sparse} \\
\hline
Fitted Error & 134.28 & 134.05 & 134.13 \\
\hline
\end{tabular*}
\end{table}

Combining Figures~\ref{figeggdir1} and \ref{figegglk}, we find that a
fly laying fewer eggs from day 4 to day 10 but making it up later by
reaching average number of egg production during the period day 10 to
day 20 is expected to live longer. Since egg production is most intense
in the early stage (day 4 to 10), this suggests a cost of early
reproduction to female Medflies. One plausible explanation is that
young Medflies are still fragile and reproduction depletes the needed
nutrition for growth.

\section{Concluding remarks}\label{sec5}

In this paper, we propose a new dimension reduction method for
longitudinal data collected over discrete,
possibly random time points. There are two key steps: the first is a
nonparametric smoothing method to borrow information from sparse
longitudinal data to estimate the inverse regression function, $E(X(t)
| Y)$; the second is the regularization needed to standardize the
longitudinal covariates.
The method is simple to implement
and effective for dimension reduction, and we establish asymptotic
theory. In particular, we achieve the optimal rate of convergence for
e.d.r. directions. Although the proposed method is
inspired by the difficulties caused by sparse longitudinal
data, the approach can also handle dense data both theoretically and
practically.

The numerical performance of the new approach is examined in a
simulation study, where we compare the estimates from dense (complete)
and sparse data. While
the results for dense data are better than those for sparse
data, the estimates for sparse data still capture the main
features of the target function. Further, these estimates are
consistent with the smoothed patterns of the estimates by \citet
{Ferre05}, and our new approach has much smaller integrated variance
with comparable integrated square bias. We also illustrate the
effectiveness of the new dimension reduction approach through the
fecundity data of Medflies with survival outcome, as only one index or
two indices are needed to summarize the longitudinal covariate
information. The high correlation between the results from
the complete and sparse data further confirms the ability
of the method in borrowing information across the entire sample for the
case of sparse data.

%


In practice, one needs to choose $k$, the dimension of the e.d.r.
space. For the Medfly data, we adopted an ad hoc approach to
subjectively select $k$ based on the fraction of variance explained by
the first few dominant eigencomponents.
For functional data, a more formal procedure to use a criterion to
measure the quality of the estimates for the e.d.r. space has been
proposed in \citet{Ferre05} as a model selection tool. Another approach
based on sequential $\chi^2$-tests was investigated in \citet{LiH:10}.
The choice of $k$ for sparse data would be an interesting topic for
future research.


\begin{appendix}\label{app}
\section{Definitions}
\begin{longlist}[1.]
\item[1.]\textit{Integrated square bias} (ISB):
\[
\operatorname{ISB} = \int\bigl(E \hat\beta(t)-\beta(t)\bigr)^2\,dt \approx\sum
_{j=1}^{N_p-1} \bigl(\bar{\hat
\beta}(t_j)-\beta(t_j) \bigr)^2(t_{j+1}-t_j),
\]
%
where $\bar{\hat\beta}(t) = \frac{1}{N_s}\sum_{i=1}^{N_s}\hat{\beta
}_i(t)$, and $N_p$ is
the number of points used to approximate\vspace*{1pt} the integral. In the
simulation study, we used $N_s=100$ and $200$ and $N_p=31$.

\item[2.]\textit{Integrated variance} (IVAR):
\[
\operatorname{IVAR} = \int E\bigl[\hat\beta(t)\bigr]^2-\bigl(E\hat\beta(t)
\bigr)^2 \,dt\approx\sum_{j=1}^{N_p-1}
\bigl(\bar{\hat\beta}{}^2(t_j)-\bigl[\bar{\hat\beta
}(t_j)\bigr]^2 \bigr) (t_{j+1}-t_j),
\]
where $\bar{\hat\beta}{}^2(t) = \frac{1}{N_s}\sum_{i=1}^{N_s}[\hat{\beta
}_i(t)]^2$.\vspace*{2pt}

\item[3.]IMSE \textit{and} MISE:
\begin{eqnarray*}
\operatorname{IMSE} &=& \int E\bigl(\hat\beta(t)-\beta(t)\bigr)^2\,dt \approx\int
\frac{1}{N_s}\sum^{N_s}_{i=1} \bigl(\hat
\beta_i(t)-\beta (t)\bigr)^2\,dt
\\
&\approx& \frac{1}{N_s}\sum^{N_s}_{i=1}
\sum_{j=1}^{N_p-1}\bigl(\hat
\beta_i(t_j)-\beta (t_j)
\bigr)^2(t_{j+1}-t_j) =\operatorname{MISE}.
\end{eqnarray*}
%
\end{longlist}

\section{Proofs}

For simplicity of notation, we let $\sum_{i,j}$ stand for
$\sum_{i=1}^n\sum_{j=1}^{N_i}$ and rewrite formula~(\ref{eq:2dsmooth}) in the main paper
as
%
\begin{equation}
\bigl(\hat{m}(\mathbf{z}),\hat{\bm{\alpha}}^T\bigr)=
\mathop {\operatorname{argmin}}_{\alpha_0,\bm{\alpha}}\sum_{i=1}^{n}
\sum_{j=1}^{N_i}\bigl\{ X_{ij}-
\alpha_0-\bm{\alpha}^T \mathbf{Z}_{ij}\bigr
\}^2K_H(\mathbf{Z}_{ij}-\mathbf{z}),\label{eq:twods_d}
\end{equation}
where $\mathbf{z}=(t,y)^T$, $\mathbf{Z}_{ij}=(T_{ij},Y_i)^T$, and
$K_H(\mathbf{u})=|H|^{-1/2}K(H^{-1/2}\mathbf{u})$, that is,
$K_H(\mathbf{Z}_{ij}-\mathbf{z})=\frac{1}{h_th_y}K_2(\frac
{T_{ij}-t}{h_t},\frac{Y_i-y}{h_y})$,
$H=  \operatorname{diag} (h_t^2, h_y^2)\sim h^2 I_{2\times2}$ [from assumption
(A.4)] and $\bm{\alpha}=(\alpha_1,\alpha_2)^T$. Here, ``$\sim$''
means ``is of the same order as.''

\begin{pf*}{Proof of Lemma~\ref{lem1}}
In order to setup the matrix-vector format of (\ref{eq:twods_d}), we
define
\begin{eqnarray*}
X&=&\bigl(\mathbf{X}_1^T,\ldots,\mathbf{X}_n^T
\bigr)^T=(X_{11},\ldots ,X_{1N_1},
\ldots,X_{n1},\ldots,X_{nN_n})^T,
\\
Z&=&\bigl(\mathbf{Z}_1^T,\ldots,\mathbf{Z}_n^T
\bigr)^T=(\mathbf{Z}_{11},\ldots ,\mathbf{Z}_{1N_1},
\ldots,\mathbf{Z}_{n1},\ldots,\mathbf{Z}_{nN_n})^T,
\end{eqnarray*}
where $\mathbf{X}_i=(X_{i1},\ldots,X_{iN_i})^T$, and
$\mathbf{Z}_i=(\mathbf{Z}_{i1},\ldots,\mathbf{Z}_{iN_i})^T$.

Thus, define
\[
Z_t=\pmatrix{ 1 & (\mathbf{z}-
\mathbf{Z}_1)
\vspace*{2pt}\cr
1 & (\mathbf{z}-\mathbf{Z}_2)
\vspace*{2pt}\cr
\vdots& \vdots
\vspace*{2pt}\cr
1 & (\mathbf{z}-\mathbf{Z}_n)}
\]
and
\[
W=\operatorname{diag} \bigl\{K_H(\mathbf{Z}_1-\mathbf{z}),
\ldots,K_H(\mathbf{Z}_n-\mathbf {z})\bigr\},
\]
the estimate $ \hat{m}(t, Y) = \hat{m}(\mathbf{z}) $ is the weighted
least square solution of (\ref{eq:twods_d})
%
\begin{equation}
\hat{m}(\mathbf{z})=\hat{\alpha}_0=e_1^T
\bigl(Z_t^T W Z_t\bigr)^{-1}Z_t^{T}
W X,\label{estimateofm_d}
\end{equation}
where $e_1=[1,0,0]^T$.

We next consider the bias and variance of $\hat m(\mathbf{z})$ in
two steps.

\textit{Step} 1: The bias of $\hat m(\mathbf{z})$ 
\[
E(\hat{\alpha}_0|\mathbf{Z}_i, i=1,\ldots,n)=e_1^T
\bigl(Z_t^T W Z_t\bigr)^{-1}Z_t^{T}
W M,
\]
where $M=(m(\mathbf{Z}_1),\ldots,m(\mathbf{Z}_n))^T$. Apply Taylor
expansion of $M$ at $\mathbf{z}$,
\[
M=Z_t \pmatrix{ m(
\mathbf{z})
\vspace*{2pt}\cr
D_m(\mathbf{z})
}+ \frac{1}{2}Q_m(
\mathbf{z})+R_m(\mathbf{z}),
\]
where $D_m$ and $Q_m$ denote the first and second derivative,
respectively, and $R_m(\mathbf{z})$ is the reminder term. Thus,
%
\begin{eqnarray}\label{eq:ExpBias_d}
E(\hat{\alpha}_0|\mathbf{Z}_i, i=1,\ldots,n)&=&e_1^T
\bigl(Z_t^T W Z_t\bigr)^{-1}Z_t^{T}
W Z_t\pmatrix{ m(
\mathbf{z})
\vspace*{2pt}\cr
D_m(\mathbf{z})}
\nonumber\\
&&{}+e_1^T\bigl(Z_t^T W
Z_t\bigr)^{-1}Z_t^T W\biggl\{
\frac{1}{2}Q_m(\mathbf{z})+R_m(\mathbf{z})\biggr
\}
\\
&=&m(\mathbf{z})+\frac{1}{2}e_1^T
\bigl(Z_t^T W Z_t\bigr)^{-1}Z_t^T
W\bigl\{Q_m(\mathbf{z})+2R_m(\mathbf{z})\bigr\}.
\nonumber
\end{eqnarray}

We first show that
%
\begin{eqnarray}\label{eq:ZWZ_d}
&&Z_t^TWZ_t\nonumber\\
&&\qquad=\pmatrix{ \displaystyle\sum_{i,j}K_H(
\mathbf{Z}_{ij}-\mathbf{z})&
\displaystyle\sum_{i,j}K_H(
\mathbf {Z}_{ij}-\mathbf{z}) (\mathbf{Z}_{ij}-
\mathbf{z})^T
\vspace*{2pt}\cr
\displaystyle\sum_{i,j}K_H(\mathbf{Z}_{ij}-
\mathbf{z}) (\mathbf{Z}_{ij}-\mathbf{z})
&
\displaystyle \sum
_{i,j}K_H(\mathbf{Z}_{ij}-\mathbf{z}) (
\mathbf{Z}_{ij}-\mathbf {z}) (\mathbf{Z}_{ij}-
\mathbf{z})^T},\nonumber\\
\end{eqnarray}
is of the order $nEN$ and find
its leading term. Since $\{(N_i, \mathbf{T}_i, Y_i), i=1,\ldots,
n\}$ and
$\{\sum_{j=1}^{N_i}K_H(\mathbf{Z}_{ij}-\mathbf{z}), i=1,\ldots,
n\}$ are i.i.d. and $N_i=N_i(n)$ and $Z_{ij}=Z_{ij}(n)$ are functions
of $n$, we will apply
classical limit theorems to the triangular array $\sum_{i,j}K_H(\mathbf
{Z}_{ij}-\mathbf{z})$ in
the first entry of (\ref{eq:ZWZ_d}).

%

The expectation of $\sum_{i,j}K_H(\mathbf{Z}_{ij}-\mathbf{z})$ is
\begin{eqnarray*}
nE\Biggl(E\Biggl(\sum_{j=1}^{N}K_H(
\mathbf{Z}_{j}-\mathbf{z})|N\Biggr)\Biggr) &=&(nEN)\cdot E
\bigl(K_H(\mathbf{Z}-\mathbf{z})\bigr)
\\
&=&(nEN)\cdot\int K(\mathbf{v})g\bigl(
\mathbf{z}+H^{1/2}\mathbf{v}\bigr)\,d\mathbf {v}
\\
&=&(nEN)\cdot\bigl(g(\mathbf{z})+o_p(1)\bigr),
\end{eqnarray*}
where the last step follows from the Taylor expansion of
$g(\mathbf{z}+H^{1/2}\mathbf{v})$ at $\mathbf{z}$.


The variance of $\sum_{i,j}K_H(\mathbf{Z}_{ij}-\mathbf{z})$ is
\begin{eqnarray*}
& & \sum_{i=1}^{n} \operatorname{var}
\Biggl(\sum_{j=1}^{N_i}K_H(
\mathbf {Z}_{ij}-\mathbf{z})\Biggr)
\\
&&\qquad=n\Biggl\{E\Biggl(\sum_{j=1}^{N}\sum
_{k=1}^{N}K_H(
\mathbf{Z}_{j}-\mathbf {z})K_H(\mathbf{Z}_{k}-
\mathbf{z})\Biggr) -\Biggl(E\Biggl(\sum_{j=1}^{N}K_H(
\mathbf{Z}_{j}-\mathbf{z})\Biggr)\Biggr)^2\Biggr\}
\\
&&\qquad=I_1-\mu_n^2/n,
\end{eqnarray*}
where
\begin{eqnarray*}
I_1&=&n\Biggl\{ E\Biggl(\sum_{j=1}^{N}
\bigl(K_H(\mathbf{Z}_{j}-\mathbf{z})\bigr)^2+
\sum_{j\neq
k}^{N}K_H(
\mathbf{Z}_{j}-\mathbf{z})K_H(\mathbf{Z}_{k}-
\mathbf{z})\Biggr)\Biggr\}
\\
&=&n\Biggl\{E\Biggl(E\Biggl(\sum_{j=1}^{N}
\bigl(K_H(\mathbf{Z}_{j}-\mathbf{z})\bigr)^2|N
\Biggr)\Biggr)\\
&&\hspace*{10pt}{}+E\Biggl(E\Biggl(\sum_{j\neq
k}^{N}K_H(
\mathbf{Z}_{j}-\mathbf{z})K_H(\mathbf{Z}_{k}-
\mathbf{z})|N\Biggr)\Biggr)\Biggr\}
\\
&=&(nEN)\cdot E\bigl(K_H(\mathbf{Z}-\mathbf{z})
\bigr)^2+\bigl(nEN(N-1)\bigr)\cdot E\bigl(K_H(\mathbf{Z}-
\mathbf{z})K_H\bigl(\mathbf{Z}'-\mathbf{z}\bigr)\bigr).
\end{eqnarray*}

Condition (A.1) implies that $\operatorname{var}(\sum_{i,j}K_H(\mathbf
{Z}_{ij}-\mathbf{z}))/(nEN)^2\rightarrow0$; hence,
\[
\frac{\sum_{i,j}K_H(\mathbf{Z}_{ij}-\mathbf{z})-E(\sum_{i,j}K_H(\mathbf
{Z}_{ij}-\mathbf{z}))}{nEN}\rightarrow0\qquad \mbox{in probability},
\]
that is,
\[
\sum_{i,j}K_H(\mathbf{Z}_{ij}-
\mathbf{z})=(nEN)\cdot\bigl(g(\mathbf{z})+o_p(1)\bigr).
\]
%

The mean of the other entries in the second row of (\ref{eq:ZWZ_d}) can be
handled similarly with
\begin{eqnarray*}
\frac{1}{nEN}E\biggl\{\sum_{i,j}K_H(
\mathbf{Z}_{ij}-\mathbf{z}) (\mathbf {Z}_{ij}-\mathbf{z})
\biggr\}&=&EK_H(\mathbf{u}-\mathbf{z}) (\mathbf{u}-\mathbf {z})
\\
&=&H\nabla g\int K_2(
\mathbf{v})\mathbf{v}\mathbf{v}^T\,d\mathbf {v}+o_p(H\cdot
\mathbf{1}),
\end{eqnarray*}
where $\mathbf{1}$ is a column vector of length 2 
with all entries equal to 1, and
\begin{eqnarray*}
&&\frac{1}{nEN}E\biggl\{\sum_{i,j}K_H(
\mathbf{Z}_{ij}-\mathbf{z}) (\mathbf {Z}_{ij}-\mathbf{z}) (
\mathbf{Z}_{ij}-\mathbf{z})^T\biggr\}
\\
&&\qquad=EK_H(\mathbf{u}-\mathbf{z}) (\mathbf{u}-\mathbf{z}) (\mathbf{u}-
\mathbf {z})^T
\\
&&\qquad=Hg(
\mathbf{z})+o_p(H).
\end{eqnarray*}
%
%
The variances of these entries in the second row of
(\ref{eq:ZWZ_d}) can be dealt with as the variance of $\sum_{i,j}K_H(\mathbf{Z}_{ij}-\mathbf{z})$, so we omit the details here and
summarize the findings for the rate
of convergence and leading terms as
\[
\frac{1}{nEN}Z_t^TWZ_t= %
\pmatrix{ g(\mathbf{z})+o_p(1) & (H\nabla g)^T+o_p
\bigl((H\cdot\mathbf{1})^T\bigr) \vspace*{2pt}
\cr
H\nabla
g+o_p(H\cdot\mathbf{1}) & Hg(\mathbf{z})+o_p(H) }
.
\]

%

We next consider the inverse $(\frac{1}{nEN}Z_t^TWZ_t)^{-1}$. To
locate the leading terms, we apply the well-known formula for a matrix
inverse in block form,
\[
\pmatrix{ a & \alpha^T \vspace*{2pt}
\cr
\alpha& B }
^{-1}= %
\pmatrix{ \displaystyle\frac{1}{a}+
\frac{1}{a^2}\alpha^T\biggl(B-\frac{\alpha\alpha
^T}{a}
\biggr)^{-1}\alpha&\displaystyle -\frac{1}{a}\alpha^T\biggl(B-
\frac{\alpha\alpha
^T}{a}\biggr)^{-1}\vspace*{4pt}
\cr
\displaystyle -\biggl(B-
\frac{\alpha\alpha^T}{a}\biggr)^{-1}\alpha\frac{1}{a}&\displaystyle \biggl(B-
\frac{\alpha
\alpha^T}{a}\biggr)^{-1} }, %
\]
where $a$ is a scalar, $\alpha$ is a column vector and $B$ is a
submatrix. Applying this formula to the matrix
$\frac{1}{nEN}Z_t^TWZ_t$, we first obtain 
\begin{eqnarray*}
&&\biggl(B-\frac{\alpha\alpha^T}{a}\biggr)^{-1}\\
&&\qquad=\biggl\{Hg(
\mathbf{z})+o_p(H)-\frac
{(H\nabla
g+o_p(H\cdot\mathbf{1}))(H\nabla g+o_p(H\cdot\mathbf
{1}))^T}{g(\mathbf{z})+o_p(1)}\biggr\}^{-1}
\\
&&\qquad=\bigl(Hg(\mathbf{z})\bigr)^{-1}+o_p
\bigl(H^{-1}\bigr).
\end{eqnarray*}

Thus,
\begin{eqnarray*}
-\biggl(B-\frac{\alpha\alpha^T}{a}\biggr)^{-1}\frac{\alpha}{a}&=&-
\bigl[\bigl(Hg(\mathbf {z})\bigr)^{-1}+o_p\bigl(H^{-1}
\bigr)\bigr]\frac{H\nabla g+o_p(H\cdot\mathbf
{1})}{g(\mathbf{z})+o_p(1)}
\\
&=& -\frac{\nabla g}{g^2(\mathbf{z})}+o_p(1),
\end{eqnarray*}
and the first entry of $(\frac{1}{nEN}Z_t^TWZ_t)^{-1}$ becomes
\begin{eqnarray*}
&&\frac{1}{a}+\frac{1}{a^2}\alpha^T\biggl(B-
\frac{\alpha\alpha
^T}{a}\biggr)^{-1}\alpha
\\
&&\qquad=\frac{1}{g(\mathbf{z})+o_p(1)}\biggl[1+\frac{1}{g(\mathbf
{z})+o_p(1)}\bigl(H\nabla
g+o_p(H\cdot\mathbf{1})\bigr)^T\biggl(
\frac{\nabla g}{g^2(\mathbf
{z})}+o_p(1)\biggr)\biggr]
\\
&&\qquad=\frac{1}{g(\mathbf{z})}+o_p(1).
\end{eqnarray*}
Therefore, we obtain
%
\begin{equation}
\label{eq:zwzinv_d}\qquad \biggl(\frac{1}{nEN}Z_t^TWZ_t
\biggr)^{-1}= %
\pmatrix{ g^{-1}(
\mathbf{z})+o_p(1)&\displaystyle-\frac{(\nabla g)^T}{g^2(\mathbf
{z})}+o_p\bigl(
\mathbf{1}^T\bigr)\vspace*{2pt}
\cr
\displaystyle-\frac{\nabla g}{g^2(\mathbf{z})}+o_p(
\mathbf{1})&\bigl(Hg(\mathbf {z})\bigr)^{-1}+o_p
\bigl(H^{-1}\bigr) } %
.
\end{equation}

Finally, we consider the rate of convergence and the
leading term of $Z_t^TWQ_m(\mathbf{z})$ in (\ref{eq:ExpBias_d}).
Since
\[
Q_m(\mathbf{z})=\bigl[(\mathbf{Z}_{11}-
\mathbf{z})^TH_m(\mathbf{z}) (\mathbf {Z}_{11}-
\mathbf{z}), \ldots,(\mathbf{Z}_{nN_n}-\mathbf{z})^TH_m(
\mathbf{z}) (\mathbf {Z}_{nN_n}-\mathbf{z})\bigr]^T,
\]
%
we have
%
%
\begin{eqnarray}\quad
\label{eq:ZWQ_d}&& \frac{1}{nEN}Z_t^TWQ_m(
\mathbf{z})
\nonumber
\\[-8pt]
\\[-8pt]
\nonumber
&&\qquad= %
\pmatrix{\displaystyle \frac{1}{nEN}\sum
_{i,j}K_H(\mathbf{Z}_{ij}-\mathbf{z}) (
\mathbf {Z}_{ij}-\mathbf{z})^TH_m(\mathbf{z}) (
\mathbf{Z}_{ij}-\mathbf{z})\vspace*{2pt}
\cr
\displaystyle \frac{1}{nEN}\sum
_{i,j}K_H(\mathbf{Z}_{ij}-
\mathbf{z}) (\mathbf {Z}_{ij}-\mathbf{z})^TH_m(
\mathbf{z}) (\mathbf{Z}_{ij}-\mathbf {z}) (\mathbf{Z}_{ij}-
\mathbf{z}) }. %
\end{eqnarray}

The second term in (\ref{eq:ZWQ_d}) is $O_p(H^{3/2}\cdot\mathbf
{1})$; hence,
\begin{eqnarray*}
&&E\bigl(\hat{m}(\mathbf{z})-m(\mathbf{z})|\mathbf{T}_i,Y_i,i=1,
\ldots,n\bigr)
\\
&&\qquad=\frac{1}{2}g^{-1}(\mathbf{z})E\biggl[\frac{1}{nEN}\sum
_{i,j}K_H(\mathbf {Z}_{ij}-
\mathbf{z}) (\mathbf{Z}_{ij}-\mathbf{z})^TH_m(
\mathbf {z}) (\mathbf{Z}_{ij}-\mathbf{z})\biggr]
\\
&&\qquad\quad{}+\biggl[-\frac{(\nabla g)^T}{g^2(\mathbf{z})}+o_p\bigl(\mathbf {1}^T
\bigr)\biggr]O_p\bigl(H^{3/2}\cdot\mathbf{1}\bigr)
\\
&&\qquad=
\frac{1}{2}g^{-1}(\mathbf{z})\int K_2(\mathbf{v})
\mathbf{v}^TH^{1/2}H_m(\mathbf{z})H^{1/2}
\mathbf {v}g\bigl(\mathbf{z}+H^{1/2}\mathbf{v}\bigr)\,d
\mathbf{v}+o_p\bigl(\operatorname{tr}(H)\bigr)
\\
&&\qquad=\frac{\sigma^2}{2}\operatorname{tr}\bigl(HH_m(\mathbf{z})\bigr)+o_p
\bigl(\operatorname{tr}(H)\bigr).
\end{eqnarray*}

\textit{Step} 2: Order of the variance (\ref{eq:var}).
From (\ref{estimateofm_d}), we know that
%
\begin{eqnarray}\label
{eq:varm_d}
&&\operatorname{var}\bigl(\hat{m}(\mathbf{z})|\mathbf{T}_i,Y_i,i=1,
\ldots,n\bigr)
\nonumber
\\[-8pt]
\\[-8pt]
\nonumber
&&\qquad= e_1^T\bigl(Z_t^T
WZ_t\bigr)^{-1}Z_t W \Sigma W
Z_t\bigl(Z_t^T W Z_t
\bigr)^{-1}e_1,
\end{eqnarray}
%
where
$\Sigma=\operatorname{var}(X|\mathbf{T}_i,Y_i,i=1,\ldots,n)=\operatorname{diag}\{\Sigma
_1,\ldots,\Sigma_n\}$,
with
$\Sigma_i=\break  \operatorname{var}(\mathbf{X}_i|\mathbf{T}_i,Y_i)=\{\sigma
_{ijk}\}$.

Let $\sum_{i,j,k}$ stand for $\sum_{i=1}^{n}\sum_{j=1}^{N_i}\sum_{k=1}^{N_i}$ and define
%
\begin{equation}
\label{eq:D_d} \biggl(\frac{1}{nEN} \biggr) Z_t^TW
\Sigma WZ_t \triangleq %
\pmatrix{ d_{11} &
d_{12} \vspace*{2pt}
\cr
d_{21} & d_{22} },
\end{equation}
where
\begin{eqnarray*}
d_{11} & =& \biggl(\frac{1}{nEN} \biggr)\sum
_{i,j,k} K_H(\mathbf {Z}_{ij}-
\mathbf{z})K_H(\mathbf{Z}_{ik}-\mathbf{z})
\sigma_{ijk},
\\
d_{12} & =& \biggl(\frac{1}{nEN} \biggr)\sum
_{i,j,k} K_H(\mathbf {Z}_{ij}-
\mathbf{z})K_H(\mathbf{Z}_{ik}-\mathbf{z}) (\mathbf
{Z}_{ik}-\mathbf{z})^T\sigma_{ijk},
\\
d_{21} & =& \biggl(\frac{1}{nEN} \biggr)\sum
_{i,j,k} K_H(\mathbf {Z}_{ij}-
\mathbf{z})K_H(\mathbf{Z}_{ik}-\mathbf{z}) (\mathbf
{Z}_{ik}-\mathbf{z})\sigma_{ijk},
\\
d_{22} & = &\biggl(\frac{1}{nEN} \biggr) \sum
_{i,j,k} K_H(\mathbf{Z}_{ij}-
\mathbf{z})K_H(\mathbf{Z}_{ik}-\mathbf{z}) (\mathbf
{Z}_{ij}-\mathbf{z}) (\mathbf{Z}_{ik}-\mathbf{z})^T
\sigma_{ijk},
\end{eqnarray*}
and $\sigma_{ijk} = \operatorname{cov}(X_{ij},X_{ik}|\mathbf
{T}_i,Y_i,i=1,\ldots,n)$.


%
%
%

The first entry of (\ref{eq:D_d}) is
\begin{eqnarray*}
d_{11}&=& E\biggl[\frac{1}{nEN}\sum_{i,j,k}K_H(
\mathbf{Z}_{ij}-\mathbf {z})K_H(\mathbf{Z}_{ik}-
\mathbf{z})\sigma_{ijk}\biggr]+o_p(1)
\\
&=&E\biggl[\frac{1}{nEN}\sum_{i,j}K_H(
\mathbf{Z}_{ij}-\mathbf{z})K_H(\mathbf {Z}_{ij}-
\mathbf{z})\Psi(T_{ij},T_{ij},Y_i)\biggr]
\\
&&{}+E\Biggl[\frac{1}{nEN}\sum_{i=1}^{n}
\sum_{j\neq
k}^{N_i}K_H(
\mathbf{Z}_{ij}-\mathbf{z})K_H(\mathbf{Z}_{ik}-
\mathbf {z})\Psi(T_{ij},T_{ik},Y_i)
\Biggr]+o_p(1)
\\
&=&I_1+I_2+o_p(1),
\end{eqnarray*}
where
\begin{eqnarray*}
I_1&=&E\bigl[K_H^2(U-\mathbf{z})\Psi(T,T,Y)
\bigr]
\\
&=&|H|^{-1/2} \int\int K_2^2(
\mathbf{u})\Psi(t+h_t t_1,t+h_t
t_1, y+h_y y_1)g\bigl(\mathbf
{z}+H^{1/2}\mathbf{u}\bigr)\,d\mathbf{u}\\
&&{}+o_p
\bigl(|H|^{1/2}\bigr)
\\
&=&|H|^{-1/2} R(K_2)\Psi(t,t,y)g(\mathbf{z})+o_p
\bigl(|H|^{-1/2}\bigr),
\end{eqnarray*}
with $R(K_2)=\int K_2^2(\mathbf{u})\,d\mathbf{u}$, and
\begin{eqnarray*}
I_2 &=&\frac{EN(N-1)}{EN}\int K_H(\mathbf{u}-
\mathbf{z})K_H(\mathbf {v}-\mathbf{z})\Psi(s_1,t_1,y_1)
g_3(s_1,t_1,y_1)\,ds_1\,dt_1\,dy_1
\\
& = &\frac{EN(N-1)}{EN h_y}
\Psi(t,t,y) g_3(t,t,y)+o_p\bigl(|H|^{-1/2}\bigr),
\end{eqnarray*}
where $\mathbf{u}=(s_1,y_1)^T$, $\mathbf{v}=(t_1,y_1)^T$, $ g_3$ is the
joint distribution of $T_{ij}$, $T_{ik}$ and $Y_i$, as defined in (A.6).

Thus,
$d_{11}=|H|^{-1/2} \Psi(t,t,y) \{ R(K_2) g(\mathbf{z})+\delta
g_3(t,t,y)\}+o_p(|H|^{-1/2})$,
where $\delta= \frac{EN(N-1)}{ENh_y}|H|^{1/2}$. The limit of $\delta$
exists when assumptions (A.1) and (A.5) hold.

Similarly, we can obtain $d_{21}$ and $d_{22}$ as follows:
\begin{eqnarray*}
d_{21}& =& O_p\bigl(|H|^{-1/2}H\bigr) +
o_p\bigl(|H|^{-1/2}H\bigr),
\\
d_{22} & =& |H|^{-1/2}H\int K_2^2(
\mathbf{u})\mathbf{u}\mathbf {u}^T\,d\mathbf{u}\Psi(t,t,y)g(
\mathbf{z})+o_p\bigl(|H|^{-1/2}H\bigr),
\end{eqnarray*}
thus, (\ref{eq:D_d}) is obtained.

Finally, we obtain (\ref{eq:var}) in the main paper
by plugging the results of (\ref{eq:zwzinv_d}) and~(\ref{eq:D_d}) into
(\ref{eq:varm_d}).
\end{pf*}

\begin{pf*}{Proof of Lemma~\ref{lem2}}
Due to space limitation, the proof is provided in the supplementary material [\citet
{JianYW:13s}] of this paper.
\end{pf*}

\begin{pf*}{Proof of Lemma~\ref{lem3}}
Similar steps as in Lemmas~\ref
{lem1} and~\ref{lem2} can be adopted and
we omit the details.
\end{pf*}

\begin{pf*}{Proof of Theorem~\ref{thmm1}}
Recall that $\hat\eta_j$ is the
estimate of the $j$th standardized e.d.r. direction, and satisfies
the equation
$\hat\Gamma^{-1/2}\hat\Gamma_e\hat\Gamma^{-1/2}\hat\eta_j=\hat\lambda
_j\hat\eta_j$.
By the definition of $\Gamma^{-1/2}$ in Section~\ref{sec2.1},
$\Gamma^{-1/2}\Gamma_e\Gamma^{-1/2}$ is a nonnegative symmetric
Hilbert--Schmidt operator, which can be interpreted as the kernel
of a linear mapping on $L_2(\mathcal{I})$. As for $\Gamma$, a
spectral decomposition of $\Gamma^{-1/2}\Gamma_e\Gamma^{-1/2}$ can
be achieved from Mercer's theorem as
\[
\Gamma^{-1/2}\Gamma_e\Gamma^{-1/2}(s,t)=\sum
_{j=1}^{k}\lambda_j\eta
_j(s)\eta_j(t),
\]
where $\lambda_1 > \lambda_2 > \cdots> \lambda_k > 0$, and
$\{\eta_j\}_{j=1}^k\cup\{\eta_j\}_{j=k+1}^\infty$ will generate a
complete orthogonal
basis of $L_2(\mathcal{I})$.

Once we estimate $\Gamma$ and $\Gamma_e$ by $\hat\Gamma$ and
$\hat\Gamma_e$, the operator
$\hat\Gamma^{-1/2}\hat\Gamma_e\hat\Gamma^{-1/2}$, which is
symmetric and Hilbert--Schmidt, has the empirical expansion:
\[
\hat\Gamma^{-1/2}\hat\Gamma_e\hat\Gamma^{-1/2}(s,t)=
\sum_{j=1}^{k}\hat \lambda_j
\hat\eta_j(s)\hat\eta_j(t).
\]

Since $\{\eta_j\}_{j=1}^\infty$ is a complete orthogonal basis,
$\hat\eta_j$ may be written as $\hat\eta_j=\mathop{\sum}_{\ell\geq
1}\hat a_{j\ell}\eta_\ell$. Let
\[
\Delta=\hat\Gamma^{-1/2}\hat\Gamma_e\hat
\Gamma^{-1/2}(s,t)-\Gamma ^{-1/2}\Gamma_e
\Gamma^{-1/2}(s,t).
\]
Following similar arguments as in Lemmas~1 and~2 of
\citet{Hall05}, we can arrive at
%
\begin{eqnarray}\label{etarate1}
\Vert \hat\eta_j-\eta_j\Vert ^2&=&D_{j2}+
\lambda_j^{-2}\Vert \Delta\Vert _{(j)}^2-
\lambda _j^{-2}\biggl(\int_{\mathcal{I}\times\mathcal{I}}\Delta
\eta_j\eta _j\biggr)^2
\nonumber
\\[-8pt]
\\[-8pt]
\nonumber
&&{}+O_p
\bigl(\Vert \Delta\Vert \Vert \Delta\Vert _{(j)}^2
\bigr),
\end{eqnarray}
where
\begin{eqnarray*}
D_{j2}&=&\sum_{\ell:\ell\neq
j}\bigl\{(
\lambda_j-\lambda_\ell)^{-2}-
\lambda_j^{-2}\bigr\}\biggl(\int_{\mathcal
{I}\times\mathcal{I}}
\Delta\eta_j\eta_\ell\biggr)^2;
\\
\Vert \Delta\Vert _{(j)}^2&\triangleq& \int
_{\mathcal{I}}\biggl\{\int_{\mathcal{I}}\Delta(s,t)
\eta_j(t)\,dt\biggr\}^2\,ds
= \sum_{\ell=1}^\infty\biggl(\int
_{\mathcal{I}\times\mathcal{I}}\Delta\eta _j\eta_\ell
\biggr)^2.
\end{eqnarray*}
As $\lambda_j$ are distinct, $D_{j2}=O_p(\Vert \Delta\Vert _{(j)}^2)$.
Thus, in order to see the asymptotic performance of $\hat\eta_j$,
it suffices to evaluate $\Vert \Delta\Vert _{(j)}$ and
$\int_{\mathcal{I}\times\mathcal{I}}\Delta\eta_j\eta_j$.

Consider $\Vert \Delta\Vert _{(j)}$. Using the sandwich
technique and the fact that $\Gamma$ and $\Gamma_e$ are continuous
operators on compact support $\mathcal{I}\times\mathcal{I}$, we
have
%
%
\begin{eqnarray}
\label{1Deltaj} \Vert \Delta\Vert _{(j)}^2&=&\bigl\Vert \hat
\Gamma^{-1/2}\hat\Gamma_e\hat\Gamma ^{-1/2}(s,t)-
\Gamma^{-1/2}\Gamma_e\Gamma^{-1/2}(s,t)\bigr\Vert
_{(j)}^2
\nonumber
\\[-8pt]
\\[-8pt]
\nonumber
&=& O_p\bigl(\Vert \hat\Gamma-\Gamma\Vert _{(j)}^2+
\Vert \hat\Gamma_e-\Gamma_e\Vert _{(j)}^2
\bigr).
\end{eqnarray}

Similar arguments as in the proof of Theorem~1 in \citet{Hall05} imply
\[
\Vert \hat\Gamma-\Gamma\Vert _{(j)}^2=O_p
\biggl(\frac{1}{nh_\phi EN}+h_\phi^4\biggr).
\]

Hereafter, for simplicity, we assume that $EX(t)=0$ so that
\[
\Gamma_e=E\bigl\{E\bigl(X(s)|Y\bigr)E\bigl(X(t)|Y\bigr)\bigr\}=E
\bigl(m(s,Y)m(t,Y)\bigr).
\]
%

From the estimation procedure,
%
%
\begin{eqnarray}
\label{Gammae} \hat\Gamma_e-\Gamma_e&=&
\frac{1}{n}\sum_{i=1}^n\hat
m(s,Y_i)\hat m(t,Y_i)-\Gamma_e
\nonumber\\
&=&\Biggl\{\frac{1}{n}\sum_{i=1}^n
\hat m(s,Y_i)\hat m(t,Y_i)-\frac{1}{n}\sum
_{i=1}^n m(s,Y_i)m(t,Y_i)
\Biggr\}
\\
&&{}+\Biggl\{\frac{1}{n}\sum_{i=1}^nm(s,Y_i)m(t,Y_i)-
\Gamma_e\Biggr\}\nonumber
\\
&\triangleq&\Delta_1+\Delta_2,\nonumber
\end{eqnarray}
where $\Delta_2=\frac{1}{n}\sum_{i=1}^nm(s,Y_i)m(t,Y_i)-\Gamma
_e=O_p(\frac{1}{\sqrt{n}})$.

Therefore,
%
\begin{equation}
\label{Deltajnorm} \Vert \hat\Gamma_e-\Gamma_e\Vert
_{(j)}^2 \leq\Vert \Delta_1\Vert
_{(j)}^2+\Vert \Delta_2\Vert
_{(j)}^2 = \Vert \Delta_1\Vert
_{(j)}^2+O_p\biggl(\frac{1}{n}\biggr).
\end{equation}
%

Since $\Delta_1=\frac{1}{n}\sum_{i=1}^n\{\hat m(s,Y_i)\hat
m(t,Y_i)-m(s,Y_i)m(t,Y_i)\}$, it suffices to consider
$\Vert \delta\Vert _{(j)}^2$, where
%
\begin{equation}
\delta=\hat m(s,y)\hat m(t,y)-m(s,y)m(t,y).\label{delta}
\end{equation}

A Taylor expansion on
$X_{ik}(T_{ik},Y_i)$ at the point $(t,y)$ leads to
%
%
\begin{eqnarray}
\label{guassexp} X_{ik}(T_{ik},Y_i)&=&X_i(t,y)+(t-T_{ik})X_i^{1,t}(t,y)
+(y-Y_i)X_i^{1,y}(t,y)
\nonumber\\
&&{}+\tfrac{1}{2}(t-T_{ik},y-Y_i)\bigtriangledown
\bigl(\bigtriangledown X_i(t_{ik},y_i)\bigr)
(t-T_{ik},y-Y_i)^T
\nonumber
\\[-8pt]
\\[-8pt]
\nonumber
&&{}+O\bigl((t-T_{ik},y-Y_i)^3
\bigr)
\\
&=&X_{ik}^{[1]}+X_{ik}^{[2]}+X_{ik}^{[3]}+X_{ik}^{[4]},\nonumber
\end{eqnarray}
where $X_i^{1,\cdot}$ is defined as the first derivative of
$X_i(t,y)$ with respect to the corresponding variable,
$\bigtriangledown$ is the gradient of $X_i(t,y)$, $t_{ik}$ and
$y_i$ are between $t$ and $T_{ik}$, and $y$ and $Y_i$,
respectively. Note that $X_{ik}^{[1]}(t,y)$ is $X_i(t,y)$,
$X_{ik}^{[2]}=(t-T_{ik})X_i^{1,t}(t,y)$,
$X_{ik}^{[3]}=(y-Y_i)X_i^{1,y}(t,y)$, and $X_{ik}^{[4]}$ is the
remaining terms in (\ref{guassexp}). Although we do not claim that
$t$ has to be close to $T_{ik}$ and $y$ to $Y_i$ in the Taylor
expansion (\ref{guassexp}), only those $\{X_{ik}\}$, whose
corresponding $(T_{ik}, Y_i)$ satisfy $|t-T_{ik}|\leq h_t$ and
$|y-Y_i|\leq h_y$, will contribute to the estimation when the
kernel weights of local linear smoother are applied. This provides
the correct order of the Taylor expansion (\ref{guassexp}), to be
elaborated below.

As defined above,
$X_{ik}^{[2]}$ and $X_{ik}^{[3]}$ have the linear term of
$t-T_{ik}$ and $y-Y_i$ whose order is $h$, and $X_{ik}^{[4]}$
contains the quadratic terms with order $h^2$. Applying the local
linear smoother in (\ref{eq:twods_d}) to (\ref{guassexp}), we
obtain
%
\begin{equation}
\hat m(t,y)=\hat m^{[1]}(t,y)+\hat m^{[2]}(t,y)+\hat
m^{[3]}(t,y)+\hat m^{[4]}(t,y),\label{msplit}
\end{equation}
where $\hat m^{[\ell]}, \ell=1,\ldots, 4$ are the corresponding
smoothers on $\{X_{ik}^{[\ell]}\}$.

Let $E'$ denote expectation conditioned on $\{(\mathbf{T}_i,Y_i),
i=1,\ldots,n\}$. Combining~(\ref{delta}) to (\ref{msplit}), we
deduce that
%
%
\begin{eqnarray}
\label{deltasplit} \delta&=&\hat m(s,y)\hat m(t,y)-E'\bigl\{\hat m(s,y)
\hat m(t,y)\bigr\}+E'\bigl\{\hat m(s,y)\hat m(t,y)\bigr
\}\nonumber\\
&&{}-m(s,y)m(t,y)
\\
&=&\bigl(\hat m^{[1]}(s,y)+\cdots+\hat m^{[4]}(s,y)\bigr)
\bigl(\hat m^{[1]}(t,y)+\cdots+\hat m^{[4]}(t,y)\bigr)\nonumber
\\
&&{} -E'\bigl\{\bigl(\hat m^{[1]}(s,y)+\cdots+\hat
m^{[4]}(s,y)\bigr)\nonumber\\
&&\hspace*{27pt}{}\times  \bigl(\hat m^{[1]}(t,y)+\cdots+\hat
m^{[4]}(t,y)\bigr)\bigr\}\nonumber
\\
&&{}+E'\bigl\{\hat m(s,y)\hat m(t,y)\bigr\}-m(s,y)m(t,y)\nonumber
\nonumber\\
&=&\delta_1+\delta_2+\delta_3+
\delta_4+\delta_5+\delta_0,\nonumber
\end{eqnarray}
where
%
%
\begin{eqnarray}
\label{deltasdef} \delta_0&=&E'\bigl\{\hat m(s,y)\hat
m(t,y)\bigr\}-m(s,y)m(t,y),\nonumber
\\
\delta_1&=&\hat m^{[1]}(s,y)\hat m^{[1]}(t,y)-E'
\bigl\{\hat m^{[1]}(s,y)\hat m^{[1]}(t,y)\bigr\},\nonumber
\\
\delta_2&=&\sum-E'\sum,\nonumber\\
\eqntext{\displaystyle\mbox{where $
\sum$ is the sum of all the linear terms of $(t-T_{ik},y-Y_i)$},}
\\
\delta_3&=&\sum-E'\sum,
\nonumber
\\[-8pt]
\\[-8pt]
 \eqntext{\displaystyle\mbox{where $
\sum$ is the sum of all the quadratic terms of $(t-T_{ik},y-Y_i)$},}
\\
\delta_4&=&\sum-E'\sum,\nonumber\\
 \eqntext{\displaystyle\mbox{where $
\sum$ is the sum of all the cubic terms of $(t-T_{ik},y-Y_i)$},}
\\
\delta_5&=&\sum-E'\sum,\nonumber\\
 \eqntext{\displaystyle\mbox{where $
\sum$ is the sum of all the quartic terms of $(t-T_{ik},y-Y_i)$}.}
\end{eqnarray}

Standard arguments in the proof of Lemma~\ref{lem1} can be applied to
$\delta_0$ to show that $\delta_0$ has the same convergence rate
as the bias of $\hat m(t,y)$, that is, $O_p(h^2)$.

Based on the results and proofs in Lemmas \ref{lem1} and~\ref
{lem2}, the
same claims as in {step} (iii) of \citet{Hall05} can be made
for $\{\delta_\ell\}$, $\ell=1,\ldots,5$. As a result, and from
(\ref{deltasplit}) and (\ref{deltasdef}), we have
%
\begin{equation}
\Vert \delta\Vert _{(j)}^2=E'\Vert
\delta_1\Vert _{(j)}^2+h^4.\label{deltajnorm}
\end{equation}
Furthermore, if we define
\[
E'\Vert \delta_1\Vert _{(j)}^2=
\int_{\mathcal{I}}\int\int_{\mathcal{I}^2}E'
\bigl\{ \delta_1(s,t_1,y)\delta_1(s,t_2,y)
\bigr\}\eta_j(t_1)\eta_j(t_2)\,dt_1\,dt_2\,ds,
\]
%
and apply the properties of two-dimensional linear smoother and the
same steps as in {step}~(v) of \citet{Hall05}, we obtain
%
\begin{equation}
E'\Vert \delta_1\Vert _{(j)}^2
\sim O_p\biggl(\frac{1}{nhEN }\biggr).\label{delta1}
\end{equation}

We have now shown
\[
\Vert \Delta_1\Vert _{(j)}^2=O_p
\biggl(\frac{1}{nhEN
}\biggr)+h^4.
\]

This and (\ref{1Deltaj})--(\ref{Deltajnorm}) imply
\[
\Vert \Delta\Vert _{(j)}^2=O_p\biggl(
\frac{1}{nhEN
}+h^4\biggr)+O_p\biggl(
\frac{1}{nh_\phi EN}+h_\phi^4\biggr).
\]

Finally, we consider
$\int_{\mathcal{I}\times\mathcal{I}}\Delta\eta_j\eta_j$ which can
be dominated by
$|\int_{\mathcal{I}\times\mathcal{I}}\delta_0\eta_j\eta_j| +
|\int_{\mathcal{I}\times\mathcal{I}}(\hat\Gamma-\Gamma)\eta_j\eta_j|$,
and hence is of order $O_p(h^2)+O_p(h_\phi^2)$. Equation
(\ref{etarate}) now follows from (\ref{etarate1}).

Next, we will discuss the optimal convergent rates under three
different sampling plans. If $EN<\infty$ (longitudinal data), it is
obvious that the optimal convergent rate $n^{-2/5}$ is achieved when $h
\sim n^{-1/5}$. When $EN\rightarrow\infty$, we assume that $ENh\sim
n^{-\tau}$, where $0 \leq\tau<1/5$ [because of assumption (A.5)].
Simple calculations show that the optimal rate $n^{-r}$ is achieved
when $h \sim n^{-(1-\tau)/4}$ and $r=-(1-\tau)/2$. Thus, the proof is complete.
\end{pf*}
\end{appendix}

\section*{Acknowledgements}
The authors would like to express
gratitude for the insightful comments of three referees, the
Associate Editor, and the Editor.

\begin{supplement}[id=suppA]
\stitle{Supplement to ``Inverse regression for longitudinal data''\\}
\slink[doi]{10.1214/13-AOS1193SUPP} 
\sdatatype{.pdf}
\sfilename{aos1193\_supp.pdf}
\sdescription{We provide additional supporting information for Section~\ref{sec2.1}, for simulation studies and for data analysis.}
\end{supplement}


%

\printaddresses

\end{document}


\begin{frontmatter}

\title{Correction: Inverse regression for longitudinal data}
\begin{center}\vspace*{12pt}
\textit{Ann. Statist.} \textbf{42} (2014) 563--591
\end{center}
\pdftitle{Correction: Inverse regression for longitudinal data}
\runtitle{Correction}

\begin{aug}
\author[A]{\fnms{Ci-Ren}~\snm{Jiang}\thanksref{T1}\ead[label=e1]{cirenjiang@stat.sinica.edu.tw}}
\and
\author[B]{\fnms{Jane-Ling}~\snm{Wang}\corref{}\thanksref{T2}\ead[label=e2]{janelwang@ucdavis.edu}}
\runauthor{C.-R. Jiang and J.-L. Wang}
\affiliation{Academia Sinica and University of California,  Davis}
\address[A]{Institute of Statistical Science\\
Academia Sinica\\
Taipei, 115\\
Taiwan\\
\printead{e1}}
\address[B]{Department of Statistics\\
University of California\\
Davis, California 95616\\
USA\\
\printead{e2}}
\end{aug}
\thankstext{T1}{Supported in part by NSC Grant NSC
101-2118-M-001-013-MY2, Taiwan.}
\thankstext{T2}{Supported in part by NSF Grants DMS-04-04630,
DMS-09-06813 and DMS-12-28369.}

%
\received{\smonth{2} \syear{2015}}


\end{frontmatter}

It has been pointed out to us by L. Forzani (Universidad Nacional del
Litoral) 
that a gap exists in Lemma~2.4 in \citet{JianYW14}. In this note
we replace it with Lemma \ref{lem2.4} and Theorem~\ref{teo2.1}, based on an additional
condition on the decay rate of $\xi_i$ and an additional assumption on
the coefficients $a_{k,i}=\langle\phi_i, \eta_k \rangle$. 
Following are the specific changes.

In the abstract, the sentence ``We develop asymptotic theory\ldots rate of
convergence.'' is replaced by ``We develop asymptotic theory for the
new procedure.''

Page 565: ``Asymptotic results for the new procedure are presented in
Section~2.3, where\ldots densely sampled longitudinal data (or functional
data).'' is replaced by ``Asymptotic results for the new procedure are
presented in Section~2.3.''

Page 579: ``In particular, we achieve the optimal rate of convergence
for e.d.r. directions'' is removed.

The text from Lemma 2.4 (page 573) to the end of page 574 is replaced
by the following:

The following condition about $\xi_i$ is similar to what has been
assumed in \citet{HallH07}.

\begin{cond}\label{cond1}
The covariance operator satisfies $\Gamma(s,t) = \sum_{i=1}^\infty\xi
_i \phi_i(s)\phi_i(t)$, where $\xi_i - \xi_{i+1} > C_0 i^{-\alpha_1-1}$
for some constant $C_0>0$, $i\geq1$ and $\alpha_1>1$; hence $\xi_i
\geq C_1 i^{-\alpha_1}$ for some constant $C_1>0$.
\end{cond}


\renewcommand{\theass}{A.9}
\begin{ass}\label{assA.9}
$\llvert a_{k,j}\rrvert < C_2 j^{-\alpha_2}$
for some constant $C_2>0$.
\end{ass}

Because $\eta_k = \sum_{i=1}^\infty a_{k,i} \phi_i(t)$ and $\eta_k \in
R_{\Gamma^{1/2}}$, $\sum_{i=1}^\infty a_{k,i}^2\xi_i^{-1} <\infty$.
Assumption~\ref{assA.9} and $\eta_k \in R_{\Gamma^{1/2}}$ imply that $2\alpha_2
- \alpha_1>1$, which means that the space spanned by $\{\eta_k\}
_{k=1}^K$ is smoother than $\Gamma$. The following two results in
\citet{HallH07} will be used in the proof.

Let $\Sigma_1 = \sum_i \tau_{1,i} \psi_{1,i}(s)\psi_{1,i}(t) $ and
$\Sigma_2 = \sum_i \tau_{2,i} \psi_{2,i}(s)\psi_{2,i}(t)$ be two
general covariance operators in $L_2$, then
%
%
\begin{eqnarray}
\label{eq:eigv} \sup_{j\geq1} \llvert\tau_{1,j} -
\tau_{2,j}\rrvert&\leq& \biggl(\int\!\!\!\int(\Sigma_1-
\Sigma_2)^2 \biggr)^{1/2};
\\
\label{eq:eigf} \sup_{j\geq1} \delta_j \llVert
\psi_{1,j} -\psi_{2,j}\rrVert&\leq& 8^{1/2} \biggl(\int
\!\!\!\int(\Sigma_1-\Sigma_2)^2
\biggr)^{1/2},
\end{eqnarray}
where $\delta_j = \min_{1\leq k\leq j} (\tau_{1,k}-\tau_{1, k+1})$ and
$\int\psi_{1,j}(t)\psi_{2,j}(t)\geq0$.

Let $\Gamma_L^{-1/2} = \sum_{i=1}^L \xi_i^{-1/2} \phi_i(s) \phi_i(t)$
be the truncated version of $\Gamma^{-1/2}$. Denote $r_{n1} =
(nh^2EN)^{-1}+ h^4$ and $r_{n2} = (nh_\mu^2EN)^{-1}+(nh_\phi^2EN)^{-1}+
(h_\mu+ h_\phi)^4$. In order to show the convergence rate of $\llVert
\hat\beta_k - \beta_k\rrVert ^2$, we need the following lemma.

\renewcommand{\thelem}{2.4}
%
\begin{lem}\label{lem2.4}
Under Condition~\ref{cond1} and Assumptions \textup{A.1}--\ref{assA.9}, we have
\begin{eqnarray*}
&& \bigl\llVert \hat\Gamma_L^{-1}\hat\Gamma_e
\hat\Gamma_L^{-1/2} - \Gamma ^{-1}
\Gamma_e\Gamma^{-1/2} \bigr\rrVert ^2
\\
&&\qquad  =
O_p \bigl(L^{(-2\alpha_2+\alpha
_1+1)} + L^{(3\alpha_1+2)} r_{n1} +
L^{(3\alpha_1-2\alpha_2+4)}r_{n2} \bigr). %
\end{eqnarray*}
\end{lem}

\renewcommand{\theteo}{2.1}
\begin{teo}\label{teo2.1}
Under Condition \ref{cond1} and Assumptions \textup{A.1}--\ref{assA.9}, we have
\[
\llVert\hat\beta_k - \beta_k\rrVert^2 =
O_p \bigl(L^{(-2\alpha_2+\alpha_1+1)} + L^{(3\alpha_1+2)} r_{n1} +
L^{(3\alpha_1-2\alpha_2+4)}r_{n2} \bigr).
\]
\end{teo}

Theorem~\ref{teo2.1} indicates that when $\alpha_2 \gg\alpha_1>1$ and all the
bandwidths are of the same order, the best convergence rates of $\hat
\beta_k$ are close to $(1/\sqrt{nh} + h^2)$ and $(1/\sqrt{nh^2} + h^2)$
for functional data (i.e., $0<EN h<\infty$) and for longitudinal data
(i.e., $EN<\infty$), respectively.


Further, the proof for Theorem~2.1 in the Appendix from page 586 to
page 590 is replaced by the following two proofs.

\begin{pf*}{Proof of Lemma \ref{lem2.4}}
Observe that
\[
\bigl\llVert\hat\Gamma_L^{-1}\hat\Gamma_e
\hat\Gamma_L^{-1/2} - \Gamma^{-1}
\Gamma_e\Gamma^{-1/2} \bigr\rrVert^2
\\
\leq3 (T_1+T_2+T_3),
\]
where\vspace*{2pt} $T_1 = \llVert \hat\Gamma_L^{-1}\hat\Gamma_e\hat\Gamma_L^{-1/2}
- \hat\Gamma_L^{-1}\Gamma_e \hat\Gamma_L^{-1/2} \rrVert ^2$\hspace*{-0.6pt},
$T_2 = \llVert \hat\Gamma_L^{-1}\Gamma_e \hat\Gamma_L^{-1/2} - \Gamma
_L^{-1}\Gamma_e\Gamma_L^{-1/2} \rrVert ^2$ and
$T_3 = \llVert \Gamma_L^{-1}\Gamma_e\Gamma_L^{-1/2} - \Gamma^{-1}\Gamma
_e\Gamma^{-1/2} \rrVert ^2$.
To complete the proof, we simply need the convergence rates of $T_1$,
$T_2$ and $T_3$. First,
\begin{eqnarray*}
T_1 
&\leq& \bigl\llVert\hat
\Gamma_L^{-1} \bigr\rrVert^2 \llVert\hat
\Gamma_e- \Gamma_e \rrVert^2 \bigl\llVert
\hat \Gamma_L^{-1/2} \bigr\rrVert^2
\\
&\leq& O_p \Biggl( \Biggl(\sum_{j=1}^L
j^{2\alpha_1} \Biggr)\times r_{n1} \times \Biggl(\sum
_{j=1}^L j^{\alpha_1} \Biggr) \Biggr) \leq
O_p \bigl( L^{(3\alpha_1+2)} r_{n1} \bigr).
\end{eqnarray*}

For the convergence rate of $T_2$, it suffices to show the order of $
\llVert (\hat\Gamma_L^{-1} -\Gamma_L^{-1})\Gamma_e \Gamma
_L^{-1/2}\rrVert ^2$, as the remainder terms are of smaller orders.
Observe that
\begin{eqnarray*}
&& \bigl\llVert \bigl(\hat\Gamma_L^{-1} -
\Gamma_L^{-1} \bigr)\Gamma_e
\Gamma_L^{-1/2} \bigr\rrVert^2
\\
&&\qquad \leq C \Biggl(
\Biggl\llVert\sum_{i=1}^L \biggl(
\frac{1}{\hat\xi_i}-\frac{1}{\xi_i} \biggr) \hat\phi_i(s) \hat
\phi_i(t) \Biggr\rrVert^2 + \Biggl\llVert\sum
_{i=1}^L \frac{1}{\xi_i} \bigl(\hat
\phi_i(s)\hat\phi_i(t) - \phi_i(s)
\phi_i(t) \bigr) \Biggr\rrVert^2 \Biggr)
\\
&&\quad\qquad{} \times \bigl\llVert\Gamma_e \Gamma_L^{-1/2}
\bigr\rrVert^2,
\end{eqnarray*}
for some positive constant $C$. Direct calculations lead to
%
%
\begin{eqnarray}
\label{eq:T2_1} \Biggl\llVert\sum_{i=1}^L
\biggl(\frac{1}{\hat\xi_i}-\frac{1}{\xi_i} \biggr) \hat\phi_i(s) \hat
\phi_i(t) \Biggr\rrVert^2 &\leq&\sum
_{i=1}^L \biggl\llvert \biggl(\frac{\xi_i -\hat\xi_i}{\hat\xi_i\xi_i}
\biggr) \biggr\rrvert^2 \leq O_p \bigl( L^{(4\alpha_1+1)}
r_{n2} \bigr),
\\
\qquad \Biggl\llVert\sum_{i=1}^L
\frac{\hat\phi_i(s) \hat\phi_i (t) - \phi_i(s) \phi_i(t) }{\xi_i} \Biggr\rrVert^2 &\leq& O_p \Biggl( \sum
_{i=1}^L \frac{\llVert \hat\phi_i - \phi_i\rrVert ^2}{\xi_i^2} \Biggr)
\nonumber\\[-8pt]\label{eq:T2_2} \\[-8pt]   \nonumber
&\leq& O_p \bigl( L^{(4\alpha_1+3)} r_{n2} \bigr)
\end{eqnarray}
and
%
%
\begin{eqnarray}\label{eq:T2_3}
\bigl\llVert\Gamma_e \Gamma_L^{-1/2}
\bigr\rrVert^2 &=& \bigl\llVert\Gamma_L^{1/2}
\Gamma_L^{-1/2}\Gamma_e \Gamma_L^{-1/2}
\bigr\rrVert^2 \leq O_p \Biggl( \sum
_{k=1}^K\sum_{i=1}^L
a_{k,i}^2 \xi_i \Biggr)
\nonumber\\[-8pt]\\[-8pt]\nonumber
&=&  O_p \bigl(L^{(-2\alpha_2 - \alpha_1+1)} \bigr).
\end{eqnarray}

Combining (\ref{eq:T2_1}), (\ref{eq:T2_2}) and (\ref{eq:T2_3}), we have
\[
 \bigl\llVert \bigl(\hat\Gamma_L^{-1} -
\Gamma_L \bigr)\Gamma_e \Gamma_L^{-1/2}
\bigr\rrVert^2 = O_p \bigl( L^{(3\alpha_1-2\alpha_2+4)}
r_{n2} \bigr).
\]

Also,
$T_3 \leq O_p (\sum_{k=1}^K \sum_{i>L} a_{k,i}^2 \xi_i^{-1} ) = O_p
(L^{(-2\alpha_2+\alpha_1+1)} )$.
The proof is thus complete.
\end{pf*}

\begin{pf*}{Proof of Theorem \ref{teo2.1}}
Observe that
\begin{eqnarray*}
\llVert\hat\beta_k - \beta_k\rrVert&= & \bigl\llVert
\hat\lambda_k^{-1} \hat\Gamma_L^{-1}
\hat\Gamma_e\hat\Gamma_L^{-1/2}\hat
\eta_k -\lambda_k^{-1} \Gamma^{-1}
\Gamma_e\Gamma^{-1/2}\eta_k \bigr\rrVert
\\
&\leq& \hat\lambda_k^{-1} \bigl\llVert\hat
\Gamma_L^{-1}\hat\Gamma_e\hat
\Gamma_L^{-1/2}\hat\eta_k - \Gamma^{-1}
\Gamma_e\Gamma^{-1/2}\hat\eta_k \bigr\rrVert
\\
&&{} + \bigl\llVert\hat\lambda_k^{-1}
\Gamma^{-1} \Gamma_e\Gamma^{-1/2}\hat
\eta_k - \lambda_k^{-1} \Gamma^{-1}
\Gamma_e \Gamma^{-1/2}\eta_k \bigr\rrVert.
\end{eqnarray*}

The convergence rate of the first term can be obtained by employing the
result of Lemma \ref{lem2.4}. 
The second term satisfies
\begin{eqnarray*}
&& \bigl\llVert\hat\lambda_k^{-1} \Gamma^{-1}
\Gamma_e\Gamma^{-1/2}\hat\eta_k -
\lambda_k^{-1} \Gamma^{-1}\Gamma_e
\Gamma^{-1/2}\eta_k \bigr\rrVert
\\
&&\qquad \leq \bigl\llVert\hat\lambda_k^{-1}
\Gamma^{-1}\Gamma_e\Gamma^{-1/2}(\hat
\eta_k -\eta_k) \bigr\rrVert+ \bigl\llVert \bigl(\hat
\lambda_k^{-1} -\lambda_k^{-1}
\bigr) \Gamma^{-1}\Gamma_e\Gamma^{-1/2}
\eta_k \bigr\rrVert,
\end{eqnarray*}
where the first component satisfies
\begin{eqnarray*}
&& \bigl\llVert\hat\lambda_k^{-1} \Gamma^{-1}
\Gamma_e\Gamma^{-1/2}(\hat\eta_k -
\eta_k) \bigr\rrVert^2
\\
&&\qquad  \leq\hat\lambda_k^{-2}
\bigl\llVert\Gamma^{-1}\Gamma_e\Gamma^{-1/2} \bigr
\rrVert^2 \bigl\llVert(\hat\eta_k -\eta_k)
\bigr\rrVert^2
\\
&&\qquad \leq O_p \Biggl( \Biggl(\sum_{k=1}^K
\sum_{j=1}^\infty a_{k,j}^2
\xi_j^{-1} \Biggr) \bigl\llVert(\hat\eta_k -
\eta_k) \bigr\rrVert^2 \Biggr)
\\
&&\qquad \leq O_p \bigl( \bigl\llVert(\hat\eta_k -
\eta_k) \bigr\rrVert^2 \bigr)
\\
&&\qquad \leq O_p \bigl(L^{(-2\alpha_2+\alpha_1+1)} + L^{(3\alpha_1+2)}
r_{n1} + L^{(3\alpha_1-2\alpha_2+4)}r_{n2} \bigr).
\end{eqnarray*}
The\vspace*{1pt} third inequality above results from the fact that $\sum_{k=1}^K \sum_{j=1}^\infty a_{k,j}^2 \xi_j^{-1} <\infty$, and the last inequality
results from (\ref{eq:eigf}). The rest of the proof follows from the
fact that
\begin{eqnarray*}
&& \bigl\llVert \bigl(\hat\lambda_k^{-1} -
\lambda_k^{-1} \bigr) \Gamma^{-1}
\Gamma_e\Gamma^{-1/2}\eta_k \bigr\rrVert^2
\\
&&\qquad  \leq \biggl(\frac{\hat\lambda_k -\lambda_k}{\hat\lambda_k\lambda
_k} \biggr)^2 \llVert
\lambda_k \beta_k\rrVert^2
\\
&&\qquad  \leq O_p \bigl( (\hat\lambda_k -
\lambda_k)^2 \bigr)
\\
&&\qquad \leq O_p \bigl(L^{(-2\alpha_2+\alpha_1+1)} + L^{(3\alpha_1+2)}
r_{n1} + L^{(3\alpha_1-2\alpha_2+4)}r_{n2} \bigr).
\end{eqnarray*}\upqed
\end{pf*}




%

\printaddresses